\documentclass[12pt]{article}

\usepackage{latexsym}
\usepackage{amsmath}
\usepackage{amssymb}
\usepackage{pstricks, pstricks-add}
\usepackage{pst-node}
\usepackage{enumerate}
\usepackage{subfigure}
\usepackage{multido}
\usepackage{calc}
\usepackage{enumerate}
\usepackage{subfigure}
\usepackage{color}
\usepackage{graphicx}
\usepackage{color}
\usepackage{titlesec}

\titleformat{\section}{\normalfont\normalsize\bfseries}{\thesection}{1em}{}
\renewcommand{\thesubsection}{\alph{subsection}}
\titleformat{\subsection}{\normalfont\normalsize\bfseries}{(\thesubsection)}{1em}{}
\renewcommand{\thesubsubsection}{\roman{subsubsection}}
\titleformat{\subsubsection}{\normalfont\normalsize\bfseries}{(\thesubsubsection)}{1em}{}

\makeatletter
\newtheorem{theorem}{Theorem}[section]
\newtheorem{lemma}[theorem]{Lemma}

\newtheorem{remark}[theorem]{Remark}
\newtheorem{cor}[theorem]{Corollary}
\newtheorem{prop}[theorem]{Proposition}
\newtheorem{defin}[theorem]{Definition}

\newtheorem{notation}[theorem]{Notation}
\newtheorem{example}[theorem]{Example}

\renewcommand{\theequation}{\arabic{equation}}

\newcommand{\Prf}{\noindent{\bf Proof}.\quad }

\newcommand{\qed}{\hfill$\Box$}

\newcommand{\eqnum}{\leavevmode\hfill\refstepcounter{equation}\textup{\tagform@{\theequation}}}

\makeatother
\title{Characterization of graphs \\ without even $F$-orientations}
\author{{M. Abreu, D. Labbate, F. Romaniello}\\
	{\small Universit\`a degli Studi della Basilicata -  Viale dell'Ateneo Lucano 10, Potenza - Italy }\\
	{\small marien.abreu@unibas.it, domenico.labbate@unibas.it, federico.romaniello@unibas.it}\\
	\and \\
	{J. Sheehan}\\
	{\small King's College -  Old Aberdeen AB24 3UE - Scotland}\\
	{\small j.sheehan@maths.abdn.ac.uk}}

\date{}

\begin{document}
\maketitle

	\begin{abstract}
		A graph $G$ is \textit{$1$--extendible} if every edge belongs to at least one $1$--factor of $G$.
		Let $G$ be a graph with a $1$--factor $F$. Then an \textit{even (odd) $F$--orientation} of $G$ is an orientation in which each $F$--alternating cycle has exactly an even (odd) number of edges directed in the same fixed direction around the cycle. If a graph $G$ admits an odd $F$--orientation for some $1$--factor $F$ then its admits an odd $F'$--orientation for all $1$--factors $F'$. Such graphs are called \textit{Pfaffian} and have been widely studied. A similar statement is not true for even $F$-orientations, of which little is known.
		The purpose of this paper is to achieve helpful results about even $F$-orientations. In particular, we examine the structure of $1$--extendible graphs $G$ which have no even $F$--orientation for a fixed $1$--factor  $F$ of $G$. 
		Such graphs contain a special subgraph in a family that we will call \textit{generalized Wagner graphs}. 
		%			and denote $\cal{W}$ which contain distinguished $1$-factors called $\cal{W}$-factors which together form a family denoted by $\cal{F}$.
		%			When a graph $G \in \cal{W}$ is such that $G$ has no even $F$-orientation for all $F \in \cal{F}$, we say that $G$ is {\em not $\cal{F}$-even}.
		%			We give a complete characterization of  in the case of graphs of connectivity at least four and $k$--regular graphs for $k \geq 3$.
		%proving that $G \in \cal{W}$ is not $\cal{F}$-even.
		We give a complete characterization of generalized Wagner graphs without even $F$-orientations in the case of connectivity at least four and in the case of $k$--regular graphs for $k \geq 3$.
		
		
	\end{abstract}

	\section{Introduction}

All graphs considered are finite and simple (without loops or multiple edges). We shall use the term multigraph when multiple edges are permitted. Most of our terminology is standard and can be found in many textbooks such as \cite{BM}, \cite{LP} and \cite{Tu2001}. We denote the fact that $H$ is a subgraph of $G$ by $H \leq G$.

Given a graph $G$ with vertex set $V(G)$ and edge set $E(G)$, we denote by $(u,v)$ an edge with end--vertices $u$ and $v$ in $G$. An {\em orientation} $\overrightarrow{G}$ of $G$ is an assignment of a direction to each edge of $G$. The opposite orientation of $\overrightarrow{G}$ is denoted by $\overleftarrow{G}$. If $\overrightarrow{G}$ is an orientation of $G$, $[u,v]_{\overrightarrow{G}}=:[u,v]$ indicates that the edge $(u,v)$ is directed from $u$ to $v$.  We say that $u$ is the {\em tail} and $v$ is the {\em head} of this edge with respect to $\overrightarrow{G}$. Sometimes we write $[u,v] \in E(\overrightarrow{G})$.

Let $G$ be a connected graph. A vertex $v \in G$ is a {\em cut vertex} of $G$ if $G-v$ is disconnected.
Let $S \subset V(G)$ if $G-S$ is disconnected, $S$ is said to be a {\em vertex cut} or a {\em separating set} of $G$.
The vertex connectivity $\kappa(G)$ of $G$, is the minimum size of a separating set.

A {\em $1$--factor} $F$ of $G$ is a regular spanning subgraph of degree one.
A graph $G$ is $1$--extendible if every edge belongs to at least one $1$--factor of $G$.
A $1$--factor $F$ of $G$ is said to {\em induce} a $1$--factor of a subgraph $H$ of $G$ if $E(H) \cap E(F)$ is a $1$--factor of $H$. Note that we will often identify $F$ with $E(F)$. A subgraph $J$ of a graph $G$ is {\em central} if $G-V(J)$ has a $1$--factor.
In \cite{CLM}, \linebreak $1$--extendible graphs are also called \emph{matching covered} and central subgraphs are called \emph{conformal}. Please note that connected $1$--extendible graphs are $2$--connected (cf. Lemma \ref{claim4.13}(ii)).

Let $F$ be a $1$--factor of $G$. Then a cycle $C$ is said to be {\em $F$--alternating} if $|E(C)|= 2|E(F) \cap E(C)|$. In particular, each $F$-alternating cycle has an even number of edges. An $F$--alternating cycle $C$ in an orientation $\overrightarrow{G}$ of $G$ is {\em evenly (oddly) oriented} if for either choice of direction of traversal around $C$, the number of edges of $C$ directed in the direction of traversal is even (odd). Since $C$ is even, this is clearly independent of the initial choice of direction around $C$.
Let $\overrightarrow{G}$ be an orientation of $G$ and $F$ be a $1$--factor of $G$. If every $F$--alternating cycle is evenly oriented then $\overrightarrow{G}$ is said to be an {\em even $F$--orientation of $G$}. On the other hand, if every $F$--alternating cycle is oddly oriented then $\overrightarrow{G}$ is said to be an {\em odd $F$--orientation of $G$}.

An $F$--orientation $\overrightarrow{G}$ of a graph $G$ is {\em Pfaffian} if it is odd.
It turns out that if $\overrightarrow{G}$ is a Pfaffian $F$--orientation then  $\overrightarrow{G}$ is a Pfaffian $F^*$--orientation for all $1$--factors $F^*$ of $G$ (cf.\cite[Theorem 8.3.2 (3)]{LP}). In this case we simply say that $G$ is {\em Pfaffian}.
It is well known that every planar graph is Pfaffian and that the smallest non--Pfaffian graph is the complete bipartite graph $K_{3,3}$.
The Petersen graph is a further example of a non--Pfaffian graph (cf. Lemma \ref{claim2.8}). 

The literature on Pfaffian graphs is extensive and the results often profound (see \cite{T} for a complete survey). In particular, the problem of characterizing Pfaffian bipartite graphs was posed  by P\'olya \cite{Po13}. Little \cite{L75} obtained the first such characterization in terms of a
family of forbidden subgraphs. Unfortunately, his characterization does not
give rise to a polynomial algorithm for determining whether a given
bipartite graph is Pfaffian, or for calculating the permanent of its
adjacency matrix when it is. Such a characterization was subsequently
obtained independently
by McCuaig \cite{McC,McC1}, and Robertson, Seymour and
Thomas \cite{RST}.
As a special case their result gives a
polynomial algorithm, and hence a good characterization, for
determining when a balanced bipartite graph $G$ with adjacency matrix $A$ is {\em det--extremal} i.e. it has $|det(A)|=per(A)$. For a structural characterization of det--extremal cubic bipartite graphs the reader may also refer to \cite{Thom86}, \cite{McC00}, \cite{McC1} and \cite{FJLS}.

The problem of characterizing Pfaffian general graphs seems much harder. Nevertheless, there have been found some very interesting connections in terms of {\em bricks} and {\em near bipartite graphs} (cf. e.g. \cite{FL}, \cite{LP}, \cite{NT}, \cite{T}, \cite{VY}).

A $1$--extendible non--bipartite graph $G$ is said to be {\em near bipartite} if there exist edges $e_1$ and $e_2$ such
that $G \backslash \{e_1 , e_2\}$ is $1$--extendible and bipartite.

The Pfaffian property which holds for odd $F$--orientations
does not hold for even $F$--orientations. Indeed, the Wagner graph $W$ is Pfaffian, so there is an odd orientation for each $1$--factor. On the other hand, it has an even $F_1$--orientation and no even $F_2$--orientation where $F_1$ and $F_2$ are chosen $1$--factors of $W$ (cf. Lemma \ref{claim2.6}).
%However, if a graph $G$ admits an even $F$ orientation for every $1$-factor $F$, we say that $G$ is $F$-even.

Since little is known about even $F$--orientations, the purpose of this paper is to achieve helpful results in this context. In particular, we examine the structure of $1$--extendible graphs $G$ which have no even $F$--orientation for a fixed $1$--factor $F$ of $G$ (cf. Theorem  \ref{theorem3.5MAIN}).
In the case of graphs of connectivity at least four and of $k$--regular graphs for $k \geq 3$  we give a characterization in (cf. Theorem  \ref{Thm2.4}.)

\section{Preliminaries}\label{prel}

In order to state our results we need some preliminary definitions and properties.

We denote by $P(u,v)$ a $uv$--path $(u:=u_0,u_1,\ldots,u_n=:v)$ and by $P(v,u)$ a $vu$--path $(v:=u_n,u_{n-1},\ldots,u_1,u_0=:u)$.
Suppose that $u$, $v$ and $w$ are distinct vertices of $G$ and that $P(u,v)$ is a $uv$--path and $Q(v,w)$ is a $vw$--path such that $V(P(u,v)) \cap V(Q(v,w))=\{v\}$. Then  $P(u,v)Q(v,w)$ denotes the $uw$--path formed by the {\em concatenation} of these paths.

\begin{defin}\label{orientation function}
	Let $\overrightarrow{G}$ be an orientation of $G$. We define a $(0,1)$--function $\omega:=\omega_{\overrightarrow{G}}$
	on the set of paths and cycles of $G$ as follows:
	
	(i) For any path $P:=P(u,v)=(u_0,\ldots,u_n)$
	$$
	\omega(P) \, := \, |\{i \, : \, [u_i,u_{i+1}] \in E(\overrightarrow{G}), 0\le i \le n-1 \}| \, (\mbox{mod} \, 2)
	$$
	
	Note that $\omega(P(u,v))\equiv \omega(P(v,u)) + n (\mbox{mod}\, 2)$;
	
	(ii) For any cycle $C=(u_1,\ldots,u_n,u_1)$
	$$
	\omega(C) \, := \, |\{i \, : \, [u_i,u_{i+1}] \in E(\overrightarrow{G}), 0\le i \le n-1 \}| \, (\mbox{mod} \, 2)
	$$
	
	where the suffixes are integers taken modulo n.
	
	We say that $\omega$ is {\em the orientation function associated with $\overrightarrow{G}$}.
	
	In other words, for each path $P$ (or cycle $C$), $\omega(P)$ (or $\omega(C)$) is the parity of the number of edges oriented consistently with  $\overrightarrow{G}$.
\end{defin}

As we have already noted, if $n$ is even then $\omega(C)$ is independent of any cyclic rotation of the vertices of $C$. This is not the case when $n$ is odd and so we have a slight abuse of notation in this case. Note also that when $n$ is even, $C$ is evenly oriented or oddly oriented if $\omega(C)=0$ or $\omega(C)=1$ respectively.

Suppose that $\overrightarrow{G}$ is an even (resp. odd) $F$--orientation of $G$ where $F$ is a fixed $1$--factor of $G$. Then the orientation function $\omega$ associated with $\overrightarrow{G}$ is said to be an {\em even $F$--function} (resp. {\em odd $F$--function}).

Observe that when $C$ is considered as a concatenation of paths, e.g.

$$C = (P_1(u_1,u_2)P_2(u_2,u_3), \ldots, P_n(u_n,u_1))$$

\noindent then

$$\omega(C) = \sum_{i=1}^n(P_i(u_i,u_{i+1})) \quad (\mbox{mod } 2).$$

\begin{defin}\label{zerosum} %{\em (Zero--sum sets)}
	Let $G$ be a graph with a $1$--factor $F$.
	Suppose that  ${\cal A} := \{ C_1,\ldots,C_k \}$ is a set of $F$--alternating cycles such that each edge
	of $G$ is contained in exactly an even number of elements of $\cal A$.
	%of $\bigcup_{i=1}^k E(C_i)$ occurs in an even number of $E(C_i)$, $i=1,\ldots,k$.
	Then ${\cal A}$ is said to be a {\em zero--sum $F$--set}.
	
	We say that the zero--sum $F$--set is respectively an {\em even $F$--set} or an {\em odd $F$--set} if $k$ is even or odd.
\end{defin}
The following lemma and its corollary are very useful to our purpose.
\begin{lemma}\label{claim2.1.1}\cite{EOAPG}
	Let $G$ be a graph with a $1$--factor $F$ and an odd zero--sum $F$--set $\cal{C}$$:=\{C_1,\ldots,C_k\}$. Suppose that $C_1,\ldots,C_{k_1}$ are oddly oriented and $C_{k_1+1},\ldots,C_k$ are evenly $F$--oriented in an orientation $\overrightarrow{G}$ of $G$. Let $k_2:=k-k_1$ and $0\le k_i \le k$ $(i=1,2)$. Then, $G$ cannot have an even $F$--orientation or an odd $F$--orientation if either $k_1$ or $k_2$ is odd, respectively.
\end{lemma}

\begin{cor}\label{claim2.1.2}\cite{EOAPG}
	Let $G$ be a graph with a $1$--factor $F$ and an odd $F$--set. Then $G$ cannot have both an odd $F$--orientation and an even $F$--orientation.
\end{cor}

\begin{figure}[h!]
	\begin{center}
		\includegraphics[scale=0.11]{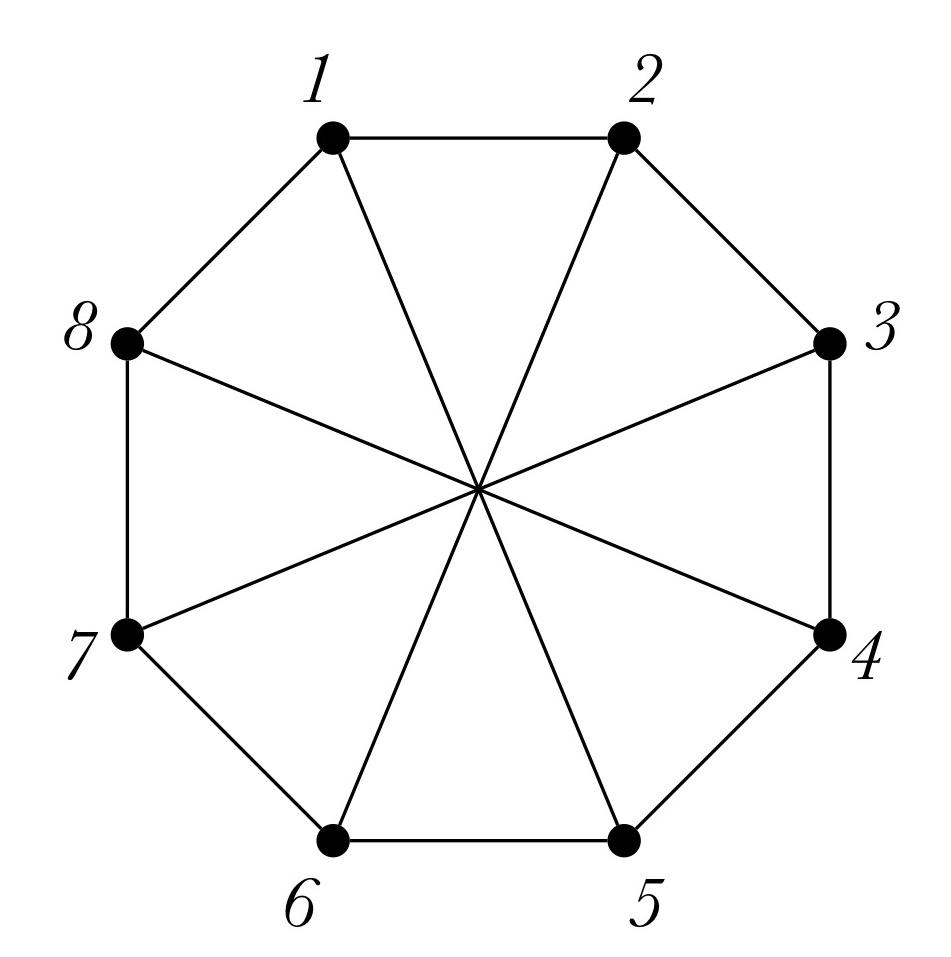}
		\caption{The \textit{Wagner Graph} $W$}
		\label{wagnerfig}	
	\end{center}
\end{figure}

The {\em Wagner graph $W$} is the cubic graph having vertex set
$V(W)=\{1,\ldots,8\}$ and edge set $E(W)$ consisting of the edges of the cycle $C=(1,\ldots,8)$ and the chords $\{(1,5),(2,6),(3,7),(4,8)\}$, see Figure \ref{wagnerfig}.

Let $C_1$ and $C_2$ be cycles of $G$ such that both include the pair of distinct independent edges $e=(u_1,u_2)$ and $f=(v_1,v_2)$. We say that $e$ and $f$ are {\em skew relative to $C_1$ and $C_2$} if the sequence $(u_1,u_2,v_1,v_2)$ occurs as a subsequence in exactly one of these cycles. Equivalently, we may write, without loss of generality, $C_1:=(u_1,u_2,\ldots,v_1,v_2,\ldots)$ and $C_2:=(u_1,u_2,\ldots,v_2,v_1,\ldots)$ i.e. if the cycles $C_1$ and $C_2$ are regarded as directed cycles, the orientation of the pair of edges $e$ and $f$ occur differently.

\begin{lemma}\label{claim2.6}
	Let $F_1:=\{(1,5)$, $(2,6)$, $(3,7)$, $(4,8)\}$ and $F_2:=\{(1,2)$, $(3,4)$, $(5,6)$, $(7,8)\}$ be $1$--factors of the Wagner graph $W$. Set $e:=(1,8)$ and $f:=(4,5)$. Then the Wagner graph $W$ satisfies the following:
	\begin{enumerate}[(i)]
		\item $W$ is $1$--extendible.
		\item $W -\{e,f\}$ is bipartite and $1$--extendible (i.e. $W$ is near bipartite).
		\item $W$ has an even $F_1$--orientation and an odd $F_1$--orientation.
		\item $W$ is Pfaffian.
		\item $W$ has no even $F_2$--orientation.
		\item There exists no pair of $F_1$--alternating cycles relative to which $e$ and $f$ are skew.
		\item The edges $e$ and $f$ are skew relative to the $F_2$--alternating cycles $C_1=(1,\ldots,8)$ and $C_2=(1,2,6,5,4,3,7,8)$.
	\end{enumerate}
\end{lemma}

\Prf
(i), (ii) and (vii) are easy to check.

(iii) The $F_1$--alternating cycles are $C_1=(1,2,6,5)$, $C_2=(2,3,7,6)$,  $C_3=(3,4,8,7)$ and $C_4=(4,5,1,8)$. It is easy to check that the orientation $\overrightarrow{W}$:
\small $$
E(\overrightarrow{W}) \, := \, \{[1,2], [2,3], [3,4], [4,5], [5,6], [6,7], [7,8], [8,1], [2,6], [3,7], [4,8], [5,1]\}
$$ \normalsize
is an even $F_1$--orientation and that the orientation $\overrightarrow{W}$
\small $$
E(\overrightarrow{W}) \, := \, \{[2,1], [2,3], [3,4], [4,5], [5,6], [6,7], [7,8], [1,8], [1,5], [2,6], [7,3], [4,8]\}
$$
is an odd $F_1$--orientation.

(iv) As we already remarked in the introduction if $W$ has an odd $F_1$--orientation then $W$ has an odd $F$--orientation for every $1$--factor $F$ of $G$. Hence, from (iii) $W$ is Pfaffian.

(v) The $F_2$--alternating cycles are:

$$C_1 = (1,2,3,4,5,6,7,8) \, , \, C_2 = (1,2,6,5,4,3,7,8)$$

$$C_3 = (1,2,3,4,8,7,6,5) \, , \, C_4 = (3,4,8,7) \, , \, C_5 = (1,2,6,5)$$

\noindent It is easy to check that $\{C_1, C_2, C_3, C_4, C_5 \}$ is an odd $F_2$--set
and that $\overrightarrow{W}$ where
$$E(\overrightarrow{W}) = \{[1,2], [2,3], [3,4], [4,5], [5,6], [6,7], [7,8], [1,8],
[5,1], [2,6], [7,3], [4,8] \}$$

is an odd $F_2$--orientation. Hence, from Corollary \ref{claim2.1.2},
%{\bf per ora e' dopo poi vedere se e come aggiustare}
$W$ has no even  $F_2$--orientation.

(vi) There is only one $F_1$--alternating cycle, namely $(4,5,1,8)$, which contains
both $e$ and $f$.
\qed

\begin{defin}\label{defin2.9}
	%(Canonical $F$--orientation)
	Let $G$ be a balanced bipartite graph with bipartition $(X,Y)$, i.e. $|X|=|Y|=n$ admitting a $1$--factor $F$. The \emph{canonical} $F$--\emph{orientation}, \noindent $\overrightarrow{G}$, of $G$ is defined as follows. Set $X:=\{x_1,x_2, \ldots, x_n\}$ and $Y:=\{y_1,y_2, \ldots, y_n\}$.
	Let $F:=\{(x_i,y_i) \, | \, i=1,2, \ldots, n \}$ be a $1$--factor of $G$.
	Let
	$$
	E(\overrightarrow{G})=\{[x_i,y_i] \, | \, i=1,2, \ldots, n\} \cup
	\{[y,x] \, | \, (y,x) \in E(G) \setminus F, x \in X, y \in Y \}
	$$
	Clearly $\overrightarrow{G}$ is an even $F$--orientation.
\end{defin}
Note that if $G$ is a bipartite graph containing a $1$--factor then $G$ has an even orientation: the canonical orientation.
In this direction, Carvalho, Lucchesi and Murty \cite{CLM2} proved that the Petersen graph and $K_{3,3}$ are non--Pfaffian.

\begin{lemma}\label{claim2.8}%\cite{CLM2}
	The Petersen graph $\cal{P}$ has an even $F$--orientation for each  $1$--factor
	$F$ of $\cal{P}$, but has no odd $F_0$--orientation, where $F_0$ is the prismatic  $1$--factor of $\cal{P}$.
\end{lemma}

\Prf
Consider the labelling of $\cal{P}$ such that  $V({\cal{P}}) := \{1, 2, \ldots, 10\}$ and  $E({\cal{P}})$ consists the cycles $(1,2,3,4,5)$, $(6,7,8,9,10)$ and the edges of the $1$--factor $F_0:=\{(1,6),(2,9),(3,7),(4,10),(5,8)\}$. It is easy to check
that the $F_0$--alternating cycles are:
$C_1=(1,6,10,4,5,8,9,2), \,  C_2=(1,6,7,3,2,9,8,5),$  $C_3=(1,6,10,4,3, \\ 7,8,5),$ $C_4=(1,6,7,3,4,10,9,2), \, C_5 = (5,8,7,3,2,9,10,4)$ 
and that the orientation $\overrightarrow{\cal{P}}:= \{[1,2],[2,3],[3,4],[4,5],[5,1],[6,7],[7,8],[8,9],[9,10],[1,6],[2,9],[3,7],$ $[4,10],[5,8],[10,6]\}$ is an even $F_0$--orientation.
Since all $1$--factors of $\cal{P}$ are similar, $\cal{P}$ has an even $F$--orientation for
all $1$--factors $F$ of $\cal{P}$. Finally, since $\{C_1,C_2,C_3,C_4,$ $C_5\}$ is an odd $F_0$--set it follows, from Corollary \ref{claim2.1.2}, that $\cal{P}$ has no odd $F_0$--orientation.
\qed

\begin{lemma}\label{claim2.11}%\cite{CLM2}
	Let $F$ be a $1$--factor of the complete bipartite graph $K_{3,3}$. Then, $K_{3,3}$ admits an even $F$--orientation but no odd $F$--orientation.
	%The complete bipartite graph $K_{3,3}$ has an even $F$--orientation but no odd $F$--orientation, where $F$ is as in Definition \eqref{defin2.9}.
\end{lemma}

\Prf
Consider $V(K_{3,3})=\{x_1,x_2,x_3,y_1,y_2,y_3\}$ where $X=\{x_1,x_2,x_3\}, Y=\{y_1,y_2,y_3\}$ \linebreak and $F=\{x_iy_i: 1 \leq i \leq 3\}$.
The $F$--alternating cycles are:
$ C_1=(x_1,y_3,x_3,y_1), \, C_2=(x_1,y_1,x_2,y_2), \, C_3=(x_2,y_3,x_3,y_2), \, C_4=(x_1,y_3,x_3,y_2,x_2,y_1),$ 
$C_5=(x_1,y_1,x_3,y_3,x_2,y_2).$
The \emph{canonical} $F$--\emph{orientation} of $K_{3,3}$ is an even $F$--orientation as seen in Definition \ref{defin2.9}.
Furthermore, it is easy to check that $\{ C_i \, | \, i=1,2, \ldots, 5\}$ is an odd
$F$--set. Hence, from Corollary \ref{claim2.1.2}, $K_{3,3}$ has no odd $F$--orientation.
\qed

\section{Main Results}\label{main}

As we have already said in the Introduction, since little is known about even $F$--orientations, the purpose of this paper is to achieve helpful results in this context.
Recall that if $G$ is a bipartite graph containing a $1$--factor then $G$ has an even orientation: the canonical orientation. We ask when graphs, not necessarily bipartite, have an even orientation. In particular, we examine the structure of $1$--extendible graphs $G$ which have no even
$F$--orientation where $F$ is a fixed $1$--factor of $G$. (cf. Theorem  \ref{theorem3.5MAIN}).

However, before stating our main theorems, again, we need some additional notation.

\begin{defin}\label{comment3.0.1}
	An {\em even subdivision} of a graph $G$ is any graph $G^*$ which can be obtained from $G$ by replacing edges $(u,v)$ of $G$ by paths $P(u,v)$ of odd length such that $V(P(u,v)) \cap V(G) = \{u,v\}$.
\end{defin}

Note that, if $F$ is a $1$--factor of $G$ then $F$ induces, in a obvious way, a $1$--factor $F^*$ of $G^*$ and conversely. For brevity, we will often blur the distinction between $F$ and $F^*$.

\begin{defin}\label{genwagner}%(Generalized Wagner graphs \cal{W})
	A graph $G$ is said to be a {\em generalized Wagner graph} if
	\begin{enumerate}[(i)]
		
		\item $G$ is $1$--extendible;
		
		\item $G$ has a subset $R:=\{e,f\}$ of edges such that $G-R$ is 1--extendible and bipartite.
		
		\item $G-R$ has a $1$--factor $F$ and there exist $F$--alternating cycles $C_1$ and $C_2$, both containing $R$, and relative to which $e$ and $f$ are skew.
	\end{enumerate}
	The set of such graphs is denoted by $\cal{W}$. We define a {\em $\cal{W}$--factor} of $G \in \cal{W}$, a $1$--factor of $G$ satisfying Definition \ref{genwagner}(iii).
	And $\cal{F}$ will denote the set of all $\cal{W}$--factors of $G \in \cal{W}$.
	
	When a graph $G \in \cal{W}$ is such that $G$ has no even $F$-orientation for all $F \in \cal{F}$, we say that $G$ is {\em not $\cal{F}$-even}.
\end{defin}

\begin{remark}\label{remark-genwagner}
	(a) For example in Lemma \ref{claim2.6},  the Wagner graph $W \in \cal{W}$ and $F_1$ is not a $\cal{W}$--factor of $W$ but $F_2$ is. Incidentally, it is easy to prove that if $G$ is a cubic graph belonging to $\cal{W}$ with at most eight vertices then $G$ is isomorphic to the Wagner graph. Thus the Wagner graph is the smallest graph in $\cal{W}$.
	
	(b) If we say that $G \in \cal{W}$ we will often assume the notation of Definition \ref{genwagner} i.e. that $F$ is a $\cal{W}$--factor of $G$ and $R$, $C_1$ and $C_2$ are as described in Definition \ref{genwagner}(ii) and (iii) respectively.
	%Often we write $\cal{W}_F$ (rather loosely) when we mean $\cal{W}_{F_0}$ where $F_0$ is a subset of $F$.
	
	(c) It is easy to see that Definition \ref{genwagner} implies that if $G \in \cal{W}$ then $G$ is near bipartite. In particular, $G$ is non--bipartite by Definition \ref{genwagner}(iii).
\end{remark}

\begin{remark}\label{evensubdivision}
	Let $G$ be a graph and $G^*$ be an even subdivision of $G$ with $F$ and $F^*$ as in Definition \ref{comment3.0.1}. Then, it is straightforward to check that $G \in \cal{W}$ if and only if $G^* \in \cal{W}$ with $F$ a $\cal{W}$--factor of $G$ and $F^*$ a $\cal{W}$--factor of $G^*$.
\end{remark}

\begin{defin}\label{w<=n}
	Let $n\ge2$ be an integer. Let $\cal{W}$$(\le n)$ denote the subset of $\cal{W}$ consisting of graphs $G$ with maximum degree $n$.
	Moreover, we define $\cal{W}$$(n)$ to be the subset of $\cal{W}$$(\le n)$ consisting of the graphs $G \in \cal{W}$$(\le n)$ such that either
	
	(i) $G$ is regular of degree $n$;
	
	or
	
	(ii) $G$ is an even subdivision of such a graph (i).
\end{defin}

Using this notation our main results are:

%%%%% MAIN THEOREM 1 %%%%%
\begin{theorem}\label{theorem3.5MAIN}
	Let $G$ be a $1$-extendible graph containing a $1$-factor $F$.
	If $G$ has no even $F$-orientation, then $G$ contains an $F$-central subgraph $H \in \cal{W}$ and $F$ is a $\cal{W}$--factor of $H$.
\end{theorem}

%%%%% MAIN THEOREM 2 - CHARACTERIZATION %%%%%
\begin{theorem}\label{Thm2.4} Let $G$ be a $1$-extendible generalized Wagner graph, i.e. $G \in \mathcal{W}$
	
	(i) If $\kappa(G) \geq 4$, then $G$ is not $\cal{F}$-even.
	
	(ii) If $k$-regular for $k \geq 3$, then $G$ is not $\cal{F}$-even.
	
	%(iii) If $G \in {\cal W}$ and it is a proper subgraph of some element of ${\cal W}^*(3)$. Then $G$ is $F$--even.
\end{theorem}

\begin{figure}[h!]
	\begin{center}
		\includegraphics[width=9cm]{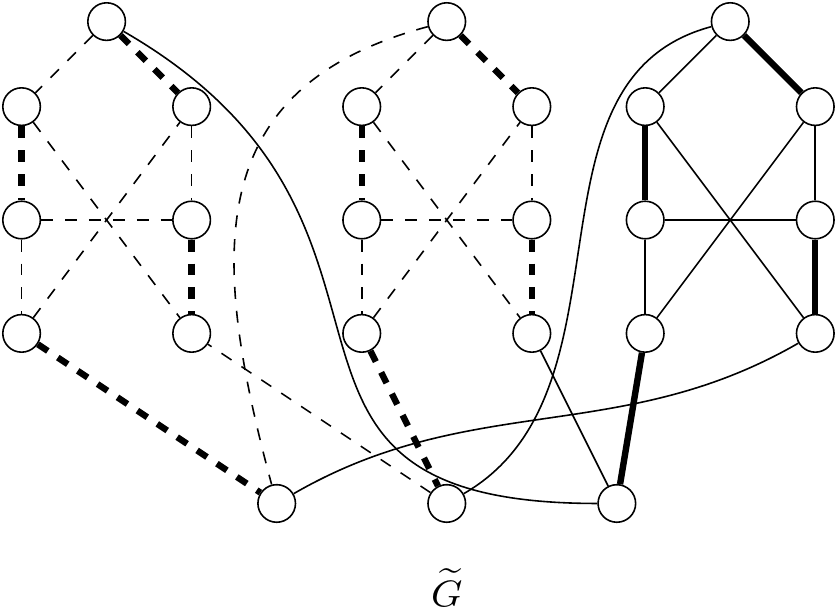}\\
	\end{center}
	\caption{$\widetilde{G} \notin \cal{W}$ but $H \in \cal{W}$ is an $F$--central subgraph}\label{exth3.5}
\end{figure}

\begin{remark}\label{commentth3.5replaced}
	Note that the graph $\widetilde{G}$ in Figure \ref{exth3.5} satisfies the conditions of Theorem \ref{theorem3.5MAIN}.
	It contains an $F$--central subgraph $H \in \cal{W}$ where the restriction of $F$ to $H$ is a $\cal{W}$--factor of $H$.
	The graph $\widetilde{G}$ admits no even $F$-orientation and $\widetilde{G} \notin \cal{W}$ because it is not near-bipartite,
	in fact at least three edges need to be removed in order to obtain a bipartite graph.
\end{remark}

The proof of these Theorems %\ref{theorem3.5MAIN}
is unfortunately very long. We begin by proving Theorem \ref{theorem3.5MAIN}. In Section \ref{1extendible} we discuss the structure of $1$--extendible graphs (see \cite{LP}). In Sections \ref{section5} and \ref{section6} the structure of a possible minimal counterexample to Theorem \ref{theorem3.5MAIN} is examined. Then in Section \ref{section7} the proof of Theorem \ref{theorem3.5MAIN} is completed.
Section \ref{section8} contains preliminary definitions and results necessary for the proof of Theorem \ref{Thm2.4} which can be found at the end of such section.

\section{Structure of 1--extendible graphs}\label{1extendible}

Let $G$ be a 1--extendible graph.
A path of odd length in $G$ whose internal vertices have degree two is called an {\em ear} of $G$. An {\em ear system} is a set $R=\{P_1,\ldots,P_n\}$ of vertex disjoint ears of $G$. Suppose that $G$ has such an ear system. Then $G-R$ is the graph obtained from $G$ by deleting all edges and the internal vertices of the constituent paths of $R$.

$R$ is said to be {\em removable} if
%\begin{enumerate}[(i)]
%\item
(i) $G-R$ is $1$--extendible
%\noindent
and
%\item
(ii) there exists no proper subset $R'$ of $R$ such that $G-R'$ is $1$--extendible.
%\end{enumerate}

\begin{defin}\label{ear decomposition}(cf. \cite{LP}, \cite{CLM})
	Let $G$ be a $1$--extendible graph. An {\em ear decomposition} of $G$ is a sequence $\cal{D}=$$(G_1, \, \ldots , \, G_r)$ of $1$--extendible graphs $G_i$ such that
	\begin{enumerate}[(i)]
		\item $G_1=K_2$, $G_r=G$;
		
		\item $G_{i-1}=G_i-R_i$, for $2\le i\le r$, where $R_i$ is a removable ear system.
	\end{enumerate}
\end{defin}

\begin{theorem}\label{claim4.2}\cite[Theorem 5.4.6]{LP}
	Let $G$ be a $1$--extendible graph and $\cal{D}=$$(G_1$, $\ldots$, $G_r)$ be an ear decomposition of $G$ with $G_{i-1}=G_i-R_i$, for $2\le i\le r$, where $R_i$ is a removable ear system. Then, for each $i$, $R_i$ has at most two ears.\qed
\end{theorem}

We say that an ear system of size $1$, size $2$ is respectively a {\em single}, {\em double} ear. If $R=\{P\}$ is a removable single ear and $P$ has length one with $E(P)=\{e\}$, then $e$ is said to be a {\em removable edge}. If $R=\{P_1, P_2\}$ is a removable double ear and $P_i$ has length one, $E(P_i)=\{e_i\}$, $i=1,2$, then $\{e_1,e_2\}$ is said to be a {\em removable doubleton}.

\begin{defin}\label{def4.3}
	Let $F$ be a $1$--factor of a $1$--extendible graph $G$. Let $\cal{D}=$$(G_1$, $\ldots$, $G_r)$ be an ear decomposition of $G$ such that $F_i:=E(F) \cap E(G_i)$ is a $1$--factor of $G_i$, $i=1,\ldots,r$. Then $\cal{D}$ is said to be an {\em $F$--reducible} ear decomposition.
\end{defin}

\begin{prop}\label{claim4.4}
	Let $F$ be a $1$--factor of a $1$--extendible graph $G$. Then there exists an $F$--reducible ear decomposition $\cal{D}=$$(G_1$, $\ldots$, $G_r)$ of $G$ with $G_{i-1}=G_i-R_i$, where $R_i$ is either a removable single ear or a removable double ear, $i=2,\ldots,r$.
\end{prop}

\Prf
We may assume that $G$ is connected. $\cal{D}$ is constructed inductively.

Let $G_1=K_2$ where $E(K_2) \subseteq E(F)$. Now suppose that for a fixed $k$, $2 \le k \le r$, there exists a sequence ${\cal{D}}_k$$=(G_1$, $\ldots$, $G_k)$ of subgraphs $G_i$ of $G$ such that, for $2\le i \le k$,
\begin{enumerate}[(i)]
	
	\item $G_{i-1}=G_i-R_i$, where $R_i$ is a removable ear system.
	
	\item $F_i$ is a $1$--factor of $G_i$ where $E(F_i)=E(F) \cap E(G_i)$.
\end{enumerate}

Suppose that $G_k \ne G$. Select, if possible, $e$ to be an edge of $G$ which has exactly one end--vertex in $G_k$. Since $G$ is $1$--extendible there exists a $1$--factor $M$ of $G$ which contains $e$. Adjoin to $G_k$ the set $R'_{k+1}$ of paths contained in $(M\backslash E(G_k)) \cup (F \backslash F_k)$. There exists at least one such path: the path containing $e$. Set $G'_{k+1}:=\bigcup R'_{k+1}$. Then $G'_{k+1}$ is $1$--extendible since $F \cap E(G'_{k+1})$ and $M \cap E(G'_{k+1})$ are both $1$--factors of $G'_{k+1}$. Now choose $R_{k+1} \subseteq R'_{k+1}$ so that $R_{k+1}$ is removable. Again $F_{k+1}:= E(F) \cap E(G_{k+1})$ is a $1$--factor of $G_{k+1}$. Thus, by induction, $\cal{D}=$$(G_1$, $\ldots$, $G_r)$ is an ear decomposition of $G$ with $G_{i-1}=G_i-R_i$, where $R_i$ is a removable ear system, for $i=2,\ldots,r$. Hence, from Theorem \ref{claim4.2}, $R_i$ has at most two ears.
Finally if $e$ cannot be chosen with exactly one end in $G_k$ then choose it so that $e$ has both ends in $G_k$, and the proof then continues exactly as in the former case.\qed

\begin{defin}\label{def4.4.1}
	%($F$--tight cut, shores)
	(i) Let $G$ be a graph and $X \subseteq V(G)$. Let $\partial(X)$ denote the set of edges with one end in $X$ and the other in $V(G) \backslash X$. A {\em cut} in $G$ is any set of the form $\partial(X)$ for some $X \subseteq V(G)$.
	
	(ii) Suppose that $G$ contains a $1$--factor $F$. A cut $\partial(X)$ is {\em $F$--tight} if $|\partial(X)\cap F|=1$. A cut is {\em tight} if it is $F$--tight for all $1$--factors $F$ of $G$. Let $G$ be a graph $G$ with a $1$--factor and $v \in V(G)$, then every cut $\partial(\{v\})$ in $G$ is tight. These tight cuts are called {\em trivial} while all the other tight cuts are called {\em non--trivial}.
	
	(iii) Let $\partial(X)$ be a non--trivial $F$--tight cut in a graph $G$  where $F$ is a $1$--factor of $G$. Let $G_1$ and $G_2$ be obtained from $G$ by identifying respectively all the vertices in $X$ and all the vertices in $\bar X :=V(G) \backslash X$ into a single vertex and deleting all resulting parallel edges. We say that $G_1$ and $G_2$ are the {\em shores} of $\partial(X)$.
	%{\bf ask john if this identification works in a single vertex}
	We denote by  $F_i$ the {\em $1$--factor of $G_i$ induced by $F$} (i=1,2).
	%(see lemma \ref{claim4.5})
\end{defin}

We now describe the Lov\'{a}sz \cite{Lo87} decomposition of $1$--extendible graphs (cf. also \cite{CLM}). Trivially we have:

\begin{lemma}\cite{Lo87},\cite{CLM}\label{claim4.5}
	Let $\partial(X)$, $X \subseteq V(G)$, be a cut in a $1$--extendible graph $G$. Then
	\begin{enumerate}[(i)]
		
		\item if $F$ is a $1$--factor of $G$, $F$ induces a $1$--factor of both of the shores of $\partial(X)$;
		
		\item if $\partial(X)$ is a tight cut then both of the shores of $\partial(X)$ are $1$--extendible.\qed
	\end{enumerate}
	
\end{lemma}

\begin{defin}\label{bricksbraces}
	%(Braces, Bricks)
	A {\em brace} (respectively a {\em brick})  is a connected bipartite (respectively a connected non--bipartite) $1$--extendible graph that has no non--trivial tight cuts.
	
	A {\em Petersen brick} is a multigraph whose undelying simple graph is the Petersen graph.
\end{defin}

\begin{defin}\label{bicritical}
	%(Bicritical)
	A graph $G$ is {\em bicritical} if $G$ contains at least one edge and $G-\{u,v\}$ has a $1$--factor for every pair of distinct vertices $u$ and $v$ in $G$.
\end{defin}

\begin{lemma}\cite{ELP82}\label{claim4.8}
	Let $G$ be a non--bipartite graph with at least four vertices. Then $G$ is a brick if and only if $G$ is $3$--connected and bicritical.\qed
\end{lemma}

%\begin{remark}\label{tightcutdecomposition}
%(Tight cut decomposition)
Let $G$ be a $1$--extendible graph with a non--trivial tight cut then, from Lemma \ref{claim4.5}, its two shores $G_1$ and $G_2$ are $1$--extendible and both are smaller than $G$. If either $G_1$ or $G_2$ has a non--trivial tight cut this procedure can be repeated. The procedure can be repeated until a list of graphs which are either bricks or braces is obtained. This is known as the {\em tight cut decomposition procedure}.
%\end{remark}

\begin{lemma}\label{claim4.10}\cite{Lo87}, \cite{CLM}
	Any two applications of the tight cut decomposition procedure yields the same list of bricks and braces, except for multiplicities of edges.\qed
\end{lemma}

\begin{lemma}\label{claim4.11}\cite{Lo87}, \cite{CLM}
	Let $G$ be a brick. If $R$ is a removable doubleton then $G-R$ is bipartite.\qed
\end{lemma}

Recall that Tutte's $1$--factor theorem states that {\em a graph $G$ has a $1$--factor if and only if $c_0(G-S) \le |S|$ for every subset $S$ of $V(G)$, where $c_0(G-S)$ denotes the number of odd components of $G-S$} (cf. e.g. \cite{BM}).
%\begin{defin}\label{def4.12}
%(Barriers)
A set $S \subseteq V(G)$ is said to be a {\em barrier} of $G$ if $c_0(G-S) > |S|$. The empty set and singletons are said to be {\em trivial} barriers.

\begin{lemma}\label{claim4.13}\cite[Theorem 1.5, Corollary 1.6]{CLM}
	
	(i) Let $G$ be a connected graph which contains a $1$--factor. Then $G$ is $1$--extendible if and only if, for every non--empty barrier $B$ of $G$, $G-B$ has no even components and no edge has both ends in $B$.
	
	(ii) Every connected $1$--extendible graph is $2$--connected. \qed
\end{lemma}

\begin{defin}\label{barriercuts-separationscuts}
	%(Barrier cuts, separation cuts)
	\begin{enumerate}[(i)]
		
		\item Suppose that $B$ is a non--trivial barrier in a connected graph $G$.
		Suppose that $H$ is a non--trivial odd component of $G-B$. Then $\partial(V(H))$ is said to be a {\em barrier cut}.
		
		\item
		Let $\{u,v\}$ $(u\ne v)$ be a non--trivial barrier, $2$--separation of a connected graph $G$. Let $G:=G_1 \cup G_2$ where $G_1 \cap G_2 = \left<u,v \right>$ (i.e. the subgraph of $G$ induced by $u$ and $v$). Then $\partial(V(G_i)-u)$, $\partial(V(G_i)-v)$ are tight cuts. Such cuts are said to be {\em 2--separation cuts} ($G-\{u,v\}$ has exactly 2 components).
	\end{enumerate}
	
\end{defin}

\begin{lemma}\label{claim4.15}\cite{ELP82}, \cite[Theorem 1.12]{CLM}
	Suppose that $G$ is a connected $1$--extendible graph which contains a non--trivial tight cut. Then $G$ has either a non--trivial barrier cut or a $2$--separation cut.\qed
\end{lemma}

\section{The structure of a minimal counterexample to Theorem \ref{theorem3.5MAIN}}\label{section5}

Let $G_0$ be such that
\begin{enumerate}[(i)]
	\item $G_0$ is a $1$--extendible graph.
	
	\item $G_0$ has no even $F$--orientation for some $1$--factor $F$ of $G_0$.
	
	\item $G_0$ contains no $F$--central subgraph $H$ such that $H \in \cal{W}$.
	
	\item $G_0$ is as small as possible subject to (i), (ii) and (iii).
\end{enumerate}
\noindent Then, if $G_0$ exists, it is a smallest counterexample to Theorem \ref{theorem3.5MAIN}.

\begin{lemma}\label{lemma5.1}
	Let $G_0$ be a smallest counterexample to Theorem \ref{theorem3.5MAIN}. Then $G_0$ is
	a non--bipartite graph and it is either $3$--connected or each $2$--separation is a barrier.
\end{lemma}
\Prf
$G_0$ is non--bipartite since otherwise $G_0$ has the canonical even $F$--orientation (see Definition \ref{defin2.9}).

By minimality $G_0$ is connected and, from Lemma \ref{claim4.13}(ii),
$G_0$ is $2$--connected.

Assume that $G_0$ has a 2--separation $\{u,v\}$ which is not a barrier. Write $G_0:=G_1 \cup G_2$ where $G_1 \cap G_2 := \left\{u,v\right\}$. Notice that, by definition, $|V(G_1)|=|V(G_2)|\equiv 0 \,$ (mod $2$), and that $G_1$ and $G_2$ are both $1$--extendible.

Let $f_1$ and $f_2$ be the edges of $F$ incident with $u$ and $v$ respectively. There are two cases to consider:

\noindent{\sc Case $(i)$:} $f_1=f_2$.

Let $F_i:=F \cap E(G_i)$. Then $F_i$ is a $1$--factor of $G_i$ ($i=1,2$). For $i=1,2$ assume that $G_i$ has an even $F_i$--orientation $\overrightarrow{G_i}$ with associated even functions $\omega_i:=\omega_{\overrightarrow{G_i}}$. We choose $\overrightarrow{G_i}$ so that $\omega_1(u,v)=\omega_2(u,v)$: this is possible since, if necessary, one can reverse all the orientations in, say, $\overrightarrow{G_1}$. Since $\{u,v\}$ is a $2$--separation, $\overrightarrow{G_1}$ and $\overrightarrow{G_2}$ together induce an even $F$--orientation of $G_0$ with associated even function $\omega_1 \cup \omega_2$. This contradicts the definition of $G_0$.

Hence, without loss of generality, we may assume that $G_1$ has no even $F_1$--orientation. By the minimality of $G_0$, $G_1$ has an $F_1$--central subgraph $H$ such that $H \in \cal{W}$. Then, it follows that $H$ is an $F$--central subgraph of $G_0$ such that $H \in \cal{W}$. Again a contradiction by the minimality of $G_0$.

\noindent{\sc Case $(ii)$:} $f_1 \ne f_2$.

Without loss of generality, we may assume that $f_1,f_2 \in E(G_1)$. Set
$$
G_i^* \, := \,
\left\{
\begin{array}{cc}
	G_i & \mbox{if \,} (u,v) \in E(G_0) \\
	G_i+(u,v) & \mbox{if \,} (u,v) \notin E(G_0)
\end{array}
\right. ,\, i=1,2 \,.
$$
Then, again, since $G_0$ is $1$--extendible and $\{u,v\}$ is a 2--separation, $G_i^*$ is $1$--extendible $(i=1,2)$.

Set $F_1:= F\cap E(G_1)$ and $F_2:=F \cap E(G_2) \cup \{(u,v)\}$. Now assume that $G_i^*$ has an even $F_i$--orientation $\overrightarrow{G_i^*}$ with associated even function $\omega_i$ $(i=1,2)$. Reversing orientations as in Case (i), if necessary, we may assume that $\omega_1(u,v)=1$ and $\omega_2(u,v)=0$.

Suppose that $C$ is any $F$--alternating cycle of $G_0$ such that $C$ is not contained in $G_i^*$ $(i=1,2)$. Then $u$ and $v$ are both vertices of $C$ since $\{u,v\}$ is a $2$--separation. Hence
$$
C \, := \, (P_1(u,v),P_2(v,u)) \,,
$$
where $P_i$ is an $F_i$--alternating path in $G_i$ $(i=1,2)$.

Again $C$ induces $F_i$--alternating cycles $C_i$ in $G_i^*$ where
$$
C_1 \, := \, (u,P_1(u,v),v)
$$
$$
C_2 \, := \, (v,P_2(v,u),u)
$$
and $\omega_i(C_i)=0$, $i=1,2$. Hence, setting $w:=\omega_1 \cup \omega_2$,
$$
\omega(C)=\omega_1(P_1(u,v))+\omega_2(P_2(v,u))=
$$
$$
=(\omega_1(P_1(u,v))+\omega_1(v,u))+(\omega_2(P_2(v,u))+\omega_2(u,v))=
$$
$$
=\omega_1(C_1)+\omega_2(C_2)=0.
$$

On the other hand, if $C$ is contained in $G_i^*$, for some $i$, then $\omega(C)=\omega_i(C)=0$ $(i=1,2)$. In all cases $\omega(C)=0$. Hence $G_0$ has an even $F$--orientation which is not the case.

Therefore, from cases (i) and (ii), we deduce that, for some $i=1,2$, $G_i^*$ has no even $F_i$--orientation.

Firstly assume that $G_i^*$ has no even $F_1$--orientation. Then, by minimality, $G_1^*$ has an $F_1$--central subgraph $H_1$ such that $H_1 \in \cal{W}$. Then, except in the case when $(u,v) \in E(H_1)$ and $(u,v) \notin E(G_0)$, $H_1$ is an $F$--central subgraph of $G_0$ such that $H_1 \in \cal{W}$. In the exceptional case, we replace $(u,v) \in E(H_1)$ by an $F_2$--alternating path $P(u,v)$ in $G_2$ to obtain an even subdivision $H_1^*$ of $H_1$ such that $H_1^*$ is an $F^*$--central subgraph of $G_0$ and $H_1^* \in \cal{W}$. Hence, using the definition of a central subgraph and Definition \ref{comment3.0.1}, again, in all cases minimality is contradicted.

Finally assume that $G_1^*$ has an even $F_1$--orientation and $G_2^*$ has no even $F_2$--orientation. The argument is almost identical as above but in the exceptional case when $(u,v)$ is, by definition, in $F_2$, and $(u,v) \notin E(G_0)$. Now as above $G_2^*$ has an $F_2$--central subgraph $H_2$ such that $H_2 \in \cal{W}$. We replace $(u,v)$ in $H_2$ by an $F_1$--alternating path in $G_1$ to obtain an even subdivision $H_2^*$ of $H_2$ such that $H_2^*$ is an $F^*$--central subgraph of $G_0$ (see Definition \ref{comment3.0.1}) and $H_2^* \in \cal{W}$. Again minimality is contradicted.

Hence, if $G_0$ is not $3$--connected each $2$--separation is a barrier.\qed

In the next lemma and subsequently, we use the notation of Definition \ref{comment3.0.1}. We need the following definition:
\begin{defin}
	Let $e_0 \in E(G)$, we say that $e \in E(G)$ is $e_0$--{\em bad} if for all $1$--factors $L$ of $G$ that contain $e$, $L$ contains $e_0$. Thus $e_0$ itself is {\em $e_0$--bad}.
\end{defin}

\begin{lemma}\label{lemma5.2}
	Let $G \in \cal{W}$ and $F$ be a $\cal{W}$--factor of $G$. Then $G$ contains an $F$--central subgraph $H$ such that $H \in \cal{W}$$(\le 3)$. Moreover $H$ is isomorphic to an even subdivision of $K_4$.
\end{lemma}

\Prf
We may assume that $G$ is connected. Suppose firstly that $G \in {\cal{W}}(3)$.
Without loss of generality $G-\{e,f\}$ is bipartite, with vertex bipartition $\{X,Y\}$ and, the edges $e$ and $f$ are skew relative to the $F$--alternating cycles $C_1$ and $C_2$. Set $e:=(x_1,x_2)$ and $f:=(y_1,y_2)$ where $x_i \in X$, $y_i \in Y$ ($i=1,2$). Set
\begin{eqnarray*}
	\nonumber %to remove numbering (before each equation)
	C_1 & = & (x_1,x_2,P_2(x_2,y_2),P_1(y_1,x_1)) \\
	C_2 & = & (x_1,x_2,Q_2(x_2,y_1),Q_1(y_2,x_1))
\end{eqnarray*}

Then we may choose $a_1, a_2 \in P_1$ and $b_1, b_2 \in P_2$
such that $Q_1(b_1,a_1)$ and $Q_2(b_2,a_2)$ are internally disjoint from $C_1$. Notice that $a_2,b_1 \in X$ and $a_1,b_2 \in Y$.
Now if $a_1 < a_2$ in $P_1(y_1,x_1)$ and $b_2 > b_1$ in $P_2(x_2,y_2)$ (or if $a_2 < a_1$ in $P_1(y_1,x_1)$ and $b_2 < b_1$ in $P_2(x_2,y_2)$) then $C_1 \cup Q_1(a_1,b_1) \cup Q_2(b_2,a_2)$ gives the required $H$. So now assume that these cases do not arise.

Hence, without loss of generality, we may assume that $a_2 < a_1$ in $P_1(y_1,x_1)$ and $b_2 > b_1$ in $P_2(x_2,y_2)$ and furthermore that $b_1$ and $b_2$ are chosen so that
\begin{itemize}
	\item $b_1 \in Q_1(y_2,x_1) \cap P_2(y_2,x_2)$ and subject to this choice $b_1$ is as large as possible in $Q_1(y_2,x_1)$ and
	\item $b_2 \in Q_2(x_2,y_1) \cap P_2(x_2,y_2)$ and subject to this choice $b_2$ is as large as possible in $Q_2(x_2,y_1)$.
\end{itemize}

Now choose $y$ in $P_1(y_1,x_1)$ so that
\begin{enumerate}[$(i)$]
	\item $y \in Q_1(y_2,x_1)$
	\item if $v > y$ in $P_1(y_1,x_1)$, $v \notin Q_2(x_2,y_1)$
	\item from $(i)$ and $(ii)$, $y$ is as small as possible in $P_1(y_1,x_1)$.
\end{enumerate}

Then choose $x \in Q_2(x_2,y_1) \cap P_1(y_1,x_1)$ so that $x < y$ in $P_1(y_1,x_1)$ and $x$ is as large as possible.

Note that by choice $x \in X$, $y \in Y$ and $P_1(x,y)$ is internally disjoint from $Q_1 \cup Q_2$. Again $P_2(b_1,b_2)$ is internally disjoint from $Q_1 \cup Q_2$. Set
$$
C_1^*:=(x_1,P_2(x_2,b_2),Q_2(b_2,y_1),P_2(y_2,b_1),Q_1(b_1,x_1)),
$$

and this case is symmetric to the one already studied with $C_1^*$,$P_1(x,y)$ and $P_2(b_1,b_2)$ taking respectively the roles of $C_1$, $Q_1(a_1,b_1)$ and $Q_2(b_2,a_2)$. Notice that now $b_2 < x$ in $Q_2(b_2,y_1)$ and $b_1 < y$ in $Q_1(b_1,x_1)$, $b_1 \in X$, $b_2 \in Y$. This gives the required $H$.

Assume now that $G$ contains a vertex $u$ with $deg(u) \ge 4$. Since $deg(u) \ge 4$ there exists $e_0:=(u,v) \in E(G)$ such that $e_0 \notin C_1\cup C_2 \cup F$. Since $G$ is $1$--extendible, $e_0 \in L_0$ for some $1$--factor $L_0$ of $G$.

Let $H$ be the graph obtained from $G$ by deleting all $e_0$--bad edges. We show that $G^* \in \cal{W}$ and $F^*$ is a $\cal{W}$--factor of $H$ (see Definition \ref{comment3.0.1} and Definition \ref{genwagner}).

\noindent{\sc Step $1$:} $C_1 \cup C_2 \subseteq H$.

Let $e \in E(C_1 \cup C_2)$. If $u \in V(C_1 \cup C_2)$ then $e$ is contained in a $1$--factor $L$ such that $e_0 \notin L$. So now suppose that $u \notin V(C_1 \cup C_2)$. If $e \in F$, then $e$ is not $e_0$--bad, since $e_0 \notin F$. Thus, w.l.o.g , we may assume that $e \in E(C_1)$ and $e \notin F$. Let $F_0$ be the $1$--factor derived from $F$ by changing the ``colours'' of $E(C_1)$.
%{\bf quest'ultima frase va aggiustata credo}
Since $u \notin V(C_1 \cup C_2)$, $e_0 \in F_0$, and $e$ is not $e_0$--bad.

\noindent{\sc Step $2$:} $H \in \cal{W}$.

Trivially $C_1$ and $C_2$ are skew relative to $e$ and $f$ in $H$ since they are skew relative to $e$ and $f$ in $G$. Furthermore, since $C_1 \cup C_2 \subseteq H$, $H - \{e,f\}$ is bipartite.

Suppose $e \in E(H)$. Then $e$ is not $e_0$--bad and hence there exists a $1$--factor $L$ of $G$ such that $e \in L$ and $e_0 \notin L$. This, in turn, implies that each edge of $L$ is not $e_0$--bad. Thus $L$ is a $1$--factor of $H$. Hence $H$ is $1$--extendible. Thus $H \in \cal{W}$, $F^*$ is a $\cal{W}$--factor of $H$ and $deg_{H}(u)=deg_G(u)-1$.

The thesis follows on repetition, if necessary, of this argument.\qed

\begin{theorem}\label{theorem5.3}
	Let $G_0$ be a minimal counterexamples to Theorem \ref{theorem3.5MAIN}. Then $G_0$ is $3$--connected.
\end{theorem}

\Prf
Assume that $G_0$ is not $3$--connected. Then, from Lemma \ref{lemma5.1}, $G_0$ has a barrier $B=\{u,v\}$, $u \ne v$. Let $H_1$ and $H_2$ be the odd components of $G_0 - B$. From Lemma \ref{claim4.13}, $G_0-B$ has no even components and $(u,v) \notin E(G_0)$. Since $G_0$ is non--bipartite at least one of $H_1$ and $H_2$ is non--trivial. So assume that $H_1$ is non--trivial  and suppose that $(u,x_1)$, $(v,y_1) \in E(F)$, $x_1 \in V(H_1)$, $y_1 \in V(H_2)$. Write $X_i:= V(H_i)$, $i=1,2$. Let $G_1$ and $G_2$ be the shores of $\partial(X_i)$ (cf. Definition \ref{def4.4.1}) where $G_1$ is obtained by contracting the vertices of $V(G_0) \backslash X_1$ to a vertex $x$ and $G_2$ is obtained contracting the vertices of $V(G_0) \backslash X_2$ to a vertex $y$.

Set $F_1:= (F \cap E(H_1)) \cup \{(x,x_1)\}$ and $F_2:= (F \cap E(H_2)) \cup \{(y,y_1)\}$. Clearly $F_i$ is a $1$--factor of $G_i$ $(i=1,2)$. From Lemma \ref{claim4.5} both $G_1$ and $G_2$ are $1$--extendible.

Since $G_0$ has no even $F$--orientation for some $i=1,2$, $G_i$ has no even $F_i$--orientation. Indeed, suppose that $G_i$ has an even orientation $\overrightarrow{G_i}$ with even $F_i$--orientation function $\omega_i$, $i=1,2$.
Set $K_1:=\partial(X)=\{(x_i,x)\, : \, i=1,\ldots, k_1 \}$ and $K_2:=\partial(Y)=\{(y_i,y)\, : \, i=1,\ldots,k_2\}$.
Moreover, suppose that $C$ is an $F$--alternating cycle of $G_0$ such that $(x_1,u)$ and $(y_1,v)$ are edges of $C$. Then $C:=(P_1(x_1,x_i),v,P_2(y_1,y_j),u)$, $2 \le i \le k_1$, $2 \le j \le k_2$, where $P_i$ is and $F_i$--alternating path in $H_i$ $(i=1,2)$.

We define an $F$--alternating function $\omega$ for $G_0$ as follows:

\noindent (i) if $(a,b) \in E(H_i)$ then $\omega(a,b)=\omega_i(a,b)$, $i=1,2$;

\noindent (ii) for edges of $E(G_0) \backslash E(H_1) \cup E(H_2)$ define

$(1) \quad \omega_1(x_i,x,x_1)+\omega_2(y_1,y,y_j):=\omega(x_1,u,y_j)+\omega(y_1,v,x_i)$.

\noindent Then, by definition of $C$, and using $(1)$:

$(2) \quad \omega(C) \, = \, \omega(P_1(x_1,x_i))+\omega(x_i,v,y_1)+\omega(P_2(y_1,y_j))+\omega(y_j,u,x_1)$

$=\omega_1((P_1(x_1,x_i))+\omega_1(x_i,x,x_1)+\omega_2(P_2(y_1,y_j))+\omega_2(y_j,y,y_1)$

$= \, \omega_1(D_1)+\omega_2(D_2)$

\noindent where $D_i$ is an $F_i$--alternating cycle in $G_i$. Hence $\omega(C) \equiv 0$ (mod $2$).

By $(i)$ if $C$ is an $F$--alternating cycle of $G_0$ not containing $(x_1,u)$ or $(y_1,v)$ then $\omega(C) \equiv 0$ (mod $2$).

This ends the proof that,  since $G_0$ has no even $F$--orientation
for some $i=1,2$, $G_i$ has no even $F_i$--orientation.

Thus, we may assume that, say $G_1$, has no even $F_1$--orientation.

By the minimality of $G_0$, $G_1$ contains an $F_1$--central subgraph $H$ such that $H \in \cal{W}$ and $F_1$ is a $\cal{W}$--factor of $H$. If $x \notin V(H)$ then $G_0$ contains $H$ and $H$ is central in $G_0$ and $F$ is a $\cal{W}$--factor of $H$, thus contradicting the minimality of $G_0$. Hence $x \in V(H)$.
By Lemma \ref{lemma5.2}, we may assume that $2\le deg_H(x) \le 3$.

Assume that $deg_H(x)=3$ and $(x,x_i) \in E(H)$, $i=1,2,3$. We may assume, without loss of generality, that either
$$
(i) \, (u,x_1),\, (u,x_2), \, (v,x_3) \, \in \, E(G_0)
$$
\noindent or
$$
(ii) \, (u,x_1),\, (v,x_2), \, (v,x_3) \, \in \, E(G_0)
$$
otherwise $H$ again would contradict the minimality of $G_0$.

We consider case (i). Let $L$ be a $1$--factor of $G_0$ containing $(v,x_3)$. Now replace the edge $(x,x_3)$ in $H$ by the path $P_1(u,x_3)$ contained in $F \cup L$ (disjoint from $H_1$) to again obtain a subgraph $H^*$ of $G_0$ with the required properties. In case (ii), Let $L$ be a $1$--factor of $G_0$ containing $(v,x_3)$. Now replace the edge $(x,x_3)$ in $H$ by the path $P_2(v,x_3)$ contained in $F \cup L$ (disjoint from $H_1$) to again obtain a subgraph $H^*$ of $G_0$ with the required properties.
Finally if $deg_H(x)=2$ then the proof of the existence of $H^*$ is exactly the same as for case (ii).

In all cases we have a contradiction with the minimality of $G_0$. Hence $G_0$ is $3$--connected.
\qed

\begin{lemma}\label{lemma5.7}
	Suppose that $G$ is a non--bipartite $1$--extendible graph with a barrier cut $B$. Let $H_1, H_2, \ldots, H_n$ ($n \ge 2$) be the odd components of $G-B$. Suppose that $G$ has no even $F$--orientation where $F$ is a $1$--factor of $G$. Set $X_i := V(H_i)$ and $G_i$ to be the shore of $\delta(X_i)$ obtained by contracting $\overline{X_i}$ to a vertex $y_i$. Set $\delta(X_i) \cap F := \{a_i,b_i\}$ where $a_i \in X_i$. Set $F_i:= (F \cap E(H_i)) \cup \{a_i,y_i\}$. Then, for some $i$, $1 \le i \le n$, $G_i$ has no even $F_i$--orientation.
\end{lemma}

\Prf
The proof follows by induction, using the argument obtained in the proof of Theorem \ref{theorem5.3}.
\qed

\begin{theorem}\label{theorem5.8}
	Let $G_0$ be a minimal counterexample to Theorem \ref{theorem3.5MAIN}. Then $G_0$ is a non--Petersen brick.
\end{theorem}

\Prf
By Lemma \ref{claim2.8} $G_0$ is not the Petersen graph. By Lemma \ref{lemma5.1} and Theorem \ref{theorem5.3}, $G_0$ is $3$--connected and not bipartite. Now suppose that $G_0$ is not a brick. Then, by definition, $G_0$ has a non--trivial tight cut. Hence, by Lemma \ref{claim4.15}, $G_0$ has a barrier cut. So by Lemma \ref{claim4.13} there exists a barrier $B$ with odd components $H_1$, $\ldots$, $H_n$ ($n\ge 2$) of $G_0-B$ such that there are no even components and $E(B)=\emptyset$. Since $G_0$ is non--bipartite, using Lemma \ref{lemma5.7} and also its notation, w.l.o.g. we assume that $H_1$ is non--trivial and that $G_1$ has no even $F_1$--orientation. Therefore, by minimality, $G_1$ has a central subgraph $H$ such that $F_1$ induced a $1$--factor and $H$ is an even subdivision of some graph in {\cal W}. As in the proof of Theorem \ref{theorem5.3}, using Lemma \ref{lemma5.2}, we may also assume that $y_1 \in V(H)$, $2\le deg_H(y_1) \le 3$ and $(y_1,a_1) \in E(H)$.

Firstly assume that $deg_H(y_1)=3$. Set $N_H:= \{x_{11},x_{12},x_{13}\}$ where $x_{11}=a_1$. Set $g_i:=(x_{1i},b_i)$ $i=1,2,3$ where $x_{11}=a_1$ and
$g_1\in F$ (recall that $F$ is a $1$--factor of $G_0$). Up to relabelling we may set $B:=\{b_1,\ldots,b_n\}$. Write $G_0^*$ for the multigraph obtained from $G_0$ by contracting each $X_i$ to a single vertex $x_i$. Clearly $G_0^*$ is a bipartite graph having the $1$--factor $F^*:=\{(x_i,b_i) | i=1,\ldots,n\}$ induced by $F$. Let $L_i$ be a $1$--factor of $G_0$ which contains $g_i$, where $L_1 \equiv F$. Notice that, since $B$ is a barrier cut, $|L_i \cap \partial(X_i)|=1$, $i=1,\ldots,n$; $j=1,2,3$.
Set $g_i^*:=(x_1,b_i)$, $i=1,2,3$. Then, $L_i$ induces naturally a $1$--factor $L_i^*$ of $G_0^*$ which contains $g_i^*$, $i=1,2,3$. Let $P_j:=P_j(b_j,b_1)$ be the $b_jb_1$--path in $L_j^* \cup F^*$ (with first edge in $F^*$), $j=2,3$. Since $b_1 \in P_2 \cap P_3$, $P_2 \cap P_3 \ne \emptyset$. Now choose $u \in V(G_0^*)$ as follows:
\begin{enumerate}[(i)]
	\item $u \in P_2 \cap P_3$;
	\item $V(P_3(b_3,u)\cap P_2)=\{b_1,u\}$, (possibly $b_1=u$).
\end{enumerate}

By construction, $u \in B$ and there exist three internally disjoint $F^*$--alternating paths $Q_j^*:=Q_j^*(u,b_j)$, $j=1,2,3$ in $G_0^*$ each of which has even length. Then, in $G_0$, we can construct three internally disjoint $F$--alternating paths $Q_j:=Q_j(u,b_j)$ from $Q_j^*$, $j=1,2,3$ as follows,
suppose that $R_j^*:=(y_1,x_i,y_2)$ is the subpath of $Q_j^*$ containing $x_i$ for some $i$, $1 \le i \le n$. We may assume that $(y_1,x_i) \in F^*$ and $(x_i,y_2) \in L_j^*$. Then there exist $x_i$, and $x_{i2}$ in $V(H_i)$ such that $(y_1, x_{i1}) \in F$ and $(x_{i2}, y_j) \in L_j$. In $Q_j^*$ we replace $R_j^*$ by the path $(y_1,R(x_i,x_{i2}),y_2)$ where $R$ is the $x_{i1}x_{i2}$--path contained in $(F \cup L_j) \cap E(H_i)$, $j=1,2$. Each of the paths $P_j^*$ if of even length. So in this way, by iteration, we obtain the required paths $Q_j(u,b_j)$, $j=1,2,3$. It follows that the graph $H_0$ defined by:
$$
V(H_0) \, = \, (V(H) \backslash \{y_1\}) \cup \{u\} \,,
$$
$$
E(H_0) \, = \, E(H-y_1) \cup Q_1 \cup Q_2 \cup Q_3 \,,
$$
is a central subgraph of $G_0$ such that $F$ induces a $1$--factor of $H_0$ and $H_0 \in \cal{W}$.

We have assumed, for the sake of clarity, that if $B^*=\{b_1,b_2,b_3\}$ then $|B^*|=3$. There is nothing to prove if $|B^*|=1$ since $H$ is already contained in $G_0$. If $|B^*|=2$ the argument is contained in the case $|B^*|=3$.

We observe that in all cases $H_0$ is contained in $G_0$ which contradicts the minimality of $G_0$. Hence $G_0$ is a non--Petersen brick.\qed

In the next theorem we use the notation of Definitions \ref{ear decomposition} and \ref{def4.3}:

\begin{theorem}\label{theorem5.10}
	Let $G_0$ be a minimal counterexample to Theorem \ref{theorem3.5MAIN}. Then $G_0$ has an $F$--reducible ear decomposition $\cal{D}=$($G_1,\ldots,G_n$),
	$(n\ge 2$; $G_0=G_n)$, such that $G_i$ has an even $F_i$--orientation $(i=1,\ldots,n-1)$ and either:
	\begin{enumerate}[(i)]
		\item $G_{n-1}=G_0 - R$, where $R=\{e\}$ is a removable edge
		%\end{enumerate}
		
		or
		
		%\begin{enumerate}[(ii)]
		\item $G_{n-1}=G_0 - R$, where $R=\{e_1,e\}$ is a removable doubleton and $G_{n-1}$ is bipartite.
	\end{enumerate}
	
\end{theorem}

\Prf
From Proposition \ref{claim4.4} $G_0$ has an $F$--reducible ear decomposition $\cal{D}=$ $(G_1$,$\ldots$, $G_n)$ with $G_n=G_0$ and $G_{i-1}=G_i-R_i$ where $R_i$ is either a removable single ear or a removable double ear. Recall that $F_i=F \cap E(G_i)$. Trivially $G_1$ ($=K_2$) has an even $F_1$--orientation. Choose $i$, $1\le i \le n$, as large as possible, so that $G_i$ has an even $F_i$--orientation. By the minimality of $G_0$, $i=n-1$. Since $G_0$ is a brick (see Theorem \ref{theorem5.8}), $G_0$ is bicritical (cf. Lemma \ref{claim4.8}). Hence, $R$ is either a removable edge or a removable doubleton. From Lemma \ref{claim4.11}, since $G_0$ is a brick, if $R$ is a removable doubleton then $G_{n-1}=G_0-R$ and $G_{n-1}$ is bipartite. \qed

\begin{remark}\label{comments5.11}
	In the next section, we prove that case $(i)$ of Theorem \ref{theorem5.10} cannot occur and we will be very close to proving Theorem \ref{theorem3.5MAIN}.
\end{remark}

\section{Theorem \ref{theorem5.10}, Case (i)}\label{section6}

We assume throughout this section that $G_0$ is a minimal counterexample to Theorem \ref{theorem3.5MAIN} and that $G_0$ has an $F$--reducible ear decomposition $\cal{D}=$ $(G_1$,$\ldots$,$G_n)$, $(n\ge2$, $G_0=G_n)$ such that $G_i$ has an even $F_i$--orientation $(i=1,\ldots,n-1)$ and $G^*:=G_{n-1}=G_0-R$ where $R=\{e\}$ is a removable edge, i.e. we assume that case (i) of Theorem \ref{theorem5.10} is true.

We now examine the structure of $G_0$ in even more detail and via a series of lemmas derive a contradiction. Our proof imitates the proof of
\cite[Theorem $1$]{L75}.

Let $\overrightarrow{G^*}$ be an even $F$--orientation of $G^*$ with associated even $F$--function $\omega$ and let $e:=(u,v)$.

\begin{lemma}\label{claim6.0}
	%Let $G^*$ be an even $F$--orientation of $G^*$ with associated even $F$--function $\omega$ and let $e:=(u,v)$.
	There exist $F$--alternating paths $Q_1:=Q_1{(u,v)}$, $Q_2:$$=Q_2{(u,v)}$ in $G^*$ such that $\omega(Q_1)$$ \ne \omega (Q_2)$. Moreover, the first and last edges of $Q_i$ $(i=1,2)$ belong to $F$.
\end{lemma}

\Prf  Since $\overrightarrow{G^*}$ is an even $F$--orientation if no such paths $Q_1$ and $Q_2$ exist, a suitable orientation of $e$ would yield an even $F$--orientation of $G_0$.

Since $e \notin F$, the first and last edges of $Q_i$ $(i=1,2)$ must belong to $F$. \qed

\begin{lemma}\label{claim6.1}
	The $F$--alternating paths $Q_1$ and $Q_2$ may be chosen in Lemma \ref{claim6.0} so that there exist $x_0$, $y_0 \in V(Q_1) \cap V(Q_2)$ such that
	\begin{enumerate}[(i)]
		\item $x_0 < y_0$ in $Q_i$ $(i=1,2)$.
		
		\item There exist paths $R_i:=R_i(x_0,y_0)$ $(i=1,2)$ such that $R_1$ and $R_2$ are respectively equal to $Q_1 \backslash Q_2$and $Q_2 \backslash Q_1$ (abusing notation slightly). The first and the last edges of $R_i$ do not belong to $F$ $(i=1,2)$.
		
		\item $\omega(R_1)=1$, $\omega(R_2)=0$;
		
		\item subject to $(i)$, $(ii)$ and $(iii)$, $|E(Q_1(u,x_0))|+|E(Q_1(y_0,v))|$ is a maximum.
		
		\item $Q_2(u,v)=Q_1(u,x_0)R_2(x_0,y_0)Q_1(y_0,v)$.
	\end{enumerate}
\end{lemma}

\Prf
Choose $Q_1$ and $Q_2$ as above and write $Q_1:=$$Q_1(a_0,$$\ldots$,$a_k)$ and $Q_2:=$$Q_2(b_0,$$\ldots$,$b_l)$, where $u=a_0=b_0$, $v=a_k=b_l$. Let $x$ be the smallest integer such that $a_x \ne b_x$. Since the first and the last edges of $Q_i$ belong to $F$, $x\ge2$ and $x \le l-2$, $x \le k-2$. Now choose $Q_1$ and $Q_2$ so that $x$ is maximized. Let $b_y$ be the first vertex of $Q_2(b_x,v)$ in $V(Q_1)$. By definition $y>x$. Set $R_1:=Q_1(a_{x-1},b_y)$, $R_2:=Q_2(a_{x-1},b_y)$, $x_0:=a_{x-1}$, $y_0:=b_y$. If $\omega(R_1) \ne \omega(R_2)$ then, without loss of generality, let $\omega(R_1)=1$ and $\omega(R_2)=0$. Finally, choose $Q_2$ such that $Q_2=Q_1(u,x_0)R_2(x_0,y_0)Q_1(y_0,v)$.

Thus we assume that $\omega(R_1)=\omega(R_2)$. Let $Q_2^*(u,v)=Q_1(u,b_y)Q_2(b_y,v)$ and replace $Q_2$ by $Q_2^*$ in the above argument. Then, by Lemma \ref{claim6.0}, the choice of $Q_1$, $Q_2$ and $x$ is contradicted.

Now choose $Q_1$, $Q_2$, $R_1$ and $R_2$ as above to maximize $|E(Q_1(u,x_0))|+|E(Q_1(y_0,v))|$. This choice implies that $Q_2(u,v)=Q_1(u,x_0)R_2(x_0,y_0)Q_1(y_0,v)$.

Note that, since $Q_1$ and $Q_2$ are $F$--alternating paths, $R_1$ and $R_2$ are $F$--alternating paths with first and last edges not in $F$.\qed

We now examine $G^*$ in more detail. Recall that $G^*=G_0-e$ and that $G^*$ is $1$--extendible.

\begin{lemma}\label{claim6.2}
	In $G^*$ there exists an edge $f$ in $R_1 \backslash F$ with the property that each $F$--alternating cycle containing $f$ has a nonempty intersection with $R_2$. Furthermore, $f$ is contained in at least one such cycle.
\end{lemma}

\Prf Suppose that the Lemma is not true. Then for each $f=(a,b) \in R_1 \backslash F$ $(a<b$ in $Q_1)$ there exists a path $P(x,y)$ $(y<a<b<x$ in $Q_1)$ where $P$ is internally disjoint from $Q_1 \cup Q_2$ and $C:=Q_1(x,y)P(x,y)$ is an $F$--alternating cycle.

Since $C$ is $F$--alternating and $Q_1$ is $F$--alternating, $Q_1(y,x)$ has first and last edge in $F$ and $P(x,y)$ has first and last edge in
$E(G^*) \backslash F$.

Let $f:=e_1=(u_1,y_0)$ where $u_1 <y_0$ in $Q_1$. From Lemma \ref{claim6.1} and the definition of $y_0$, $e_1 \in R_1 \backslash F$. Choose a path $P_1(x_1,y_1)$, $y_1<u_1<y_0<x_1$ in $Q_1$ where $P_1$ is internally disjoint from $Q_1 \cup Q_2$ and $C_1:=Q_1(y_1,x_1)P_1(x_1,y_1)$ is an $F$--alternating cycle in $G^*$. We choose $x_1$ and $y_1$ to minimize the length of $Q_1(u_1,y_1)$.

If $y_1 \in V(R_1)$, we repeat the procedure  with $y_1$ playing the role of $y_0$. In the same way we choose $y_2$, $x_2$, $P_2(x_2,y_2)$ and $C_2:=Q_1(y_2,x_2)P_2(x_2,y_2)$ such that the length of $Q_1(u_2,y)$ is minimized. Because of the minimization of the lengths of $Q_1(u_i,y_i)$, $i=1,2$:
\begin{enumerate}[(i)]
	\item $y_2<y_1<x_2<y_0<x_1$ in $Q_1$;
	\item $P_1(x_1,y_1)$ and $P_2(x_2,y_2)$ are disjoint.
\end{enumerate}

We repeat this argument and continue to construct disjoint paths $P_i:=P_i(x_i,y_i)$ and $F$--alternating cycles $C_i:=Q_1(y_i,x_i)P_i(x_i,y_i)$, ($y_{i-1}<y_{i-2}<x_{i-1}<y_{i-3} < \ldots < x_2<y_0<x_1$) until we reach an integer $j$ such that $y_j \in Q(u,x_0)$ and $y_{j-1} \in R_1(x_0,y_0)$. Since $C_j$ is $F$--alternating and the first and last edges of $P_j$ do not belong to $F$, $y_j \ne x_0$.

Now let $\overrightarrow{G^*}$ be a fixed even $F$--orientation of $G^*$ with associated even function $\omega$. Since $\omega$ is even and $C_i$ is an $F$--alternating cycle in $G^*$, $\omega(C_i)=0$, for $i=1,\ldots,j$. Hence
\setcounter{equation}{0}

\begin{equation}\label{eq(1)}
	\sum_{i=1}^j \omega(Q_1(y_i,x_i)) + \sum_{i=1}^j \omega(P_i) \equiv 0 \, (\mbox{mod} \, 2) \,.
\end{equation}
Set
\begin{equation}\label{eq(2)}
	\begin{split}
		C:=Q_1(y_j,x_0)R_2(x_0,y_0)Q_1(y_0,x_1)P_1(x_1,y_1)Q_1(y_1,x_2)P_2(x_2,y_2) & \\
		Q_1(y_2,x_3)P_3(x_3,y_3)\ldots P_{j-1}(x_{j-1},y_{j-1})Q_1(y_{j-1}x_j)P_j(x_j,y_j) \,.
	\end{split}
\end{equation}
By definition, $C$ is an $F$--alternating cycle in $G^*$ and therefore $\omega(C)=0$.
Hence, using Lemma \ref{claim6.1}(iii) and (\ref{eq(2)})
\begin{equation}\label{eq(3)}
	\omega(Q_1(y_j,x_0))+\sum_{i=1}^j\omega(Q_1(y_{j-1},x_i))+\sum_{i=1}^j \omega(P_i) \equiv 0 \, (\mbox{mod}
	\, 2) \,.
\end{equation}
Since
\begin{equation}\label{eq(4)}
	\begin{split}
		Q_1(y_i,x_i)=Q_1(y_i,y_{i-1})+Q_1(y_{i-1},x_i),\\
		\omega(Q_1(y_i,x_i)) \equiv \omega(Q_1(y_i,y_{i-1}))+ \omega(Q_1(y_{i-1},x_i)) \, \, (\mbox{mod} \, 2) \,.
	\end{split}
\end{equation}
Adding (\ref{eq(1)}) and (\ref{eq(3)})
\begin{equation}\label{eq(5)}
	\begin{split}
		\sum_{i=2}^{j-1} (\omega(Q_1(y_i,x_i)) + \omega(Q_1(y_{i-1},x_i)))+(\omega(Q_1(y_0,x_1))+\omega(Q_1(y_1,x_1))+\\
		+\omega(Q_1(y_j,x_j))+\omega(Q_1(y_j,x_0))+\omega(Q_1(y_{j-1},x_j))) \equiv 0 \, (\mbox{mod} \, 2) \,.
	\end{split}
\end{equation}
From (\ref{eq(5)}), using (\ref{eq(4)})
\begin{equation}\label{eq(6)}
	\begin{split}
		\sum_{i=2}^{j-1} (\omega(Q_1(y_i,y_{i-1})) + \omega(Q_1(y_1,x_1))+(\omega(Q_1(y_0,x_1))+\omega(Q_1(y_j,x_j))+\\
		+\omega(Q_1(y_j,x_0))+\omega(Q_1(y_{j-1},x_j))) \equiv 0 \, (\mbox{mod} \, 2) \,.
	\end{split}
\end{equation}
i.e.
$$
\omega(Q_1(y_{j-1},y_1)) + \omega(Q_1(y_1,y_0))+\omega(Q_1(x_0,y_{j-1})) \equiv 0 \, (\mbox{mod} \, 2) \,.
$$
i.e. $\omega(R_1)=0$ which contradicts Lemma \ref{claim6.1}(iii).\qed

\begin{lemma}\label{claim6.3}
	Case (i) of Theorem \ref{theorem5.10} is not possible.
\end{lemma}

\Prf
The result is proved by contradiction. Using Lemma \ref{claim6.2} we can select an edge $f:=(a,b)$ in $R_1 \backslash F$ and an $F$--alternating cycle $C$ such that for some $z,x_1 \in V(Q_1)$, $z<a<b<x_1$ $(x_1 \ne y_0)$ and $C:=Q_1(z,x_1)P(x_1,z)$ where $P(x_1,z) \cap R_2(x_0,y_0) \ne \emptyset$.

Now choose $y_1 \in V(R_2)$ $(y_1 \ne y_0)$ so that $P_1:=P(x_1,y_1)$ is edge--disjoint from $R_2$. Furthermore, choose $x_1$ and $y_1$ to minimize the length of $Q_2(u,y_1)$.

We repeat the argument of Lemma \ref{claim6.2}. In that Lemma we begin with the edge $e_1=(u_1,y_0)$ where $u_1<y_0$ in $Q_1$. We now start with the edge $e_2:=(u_1^*,y_1)$ in $Q_2$ where $e_2 \in R_2 \backslash F$. The edge $e_2$ plays the role of $e_1$ below.

As in Lemma \ref{claim6.2} we construct disjoint $F$--alternating paths $P_i:=P_i(x_i,y_i)$, $i=1,\ldots,j$ such that
\begin{enumerate}[(i)]
	\item $P_i$ is edge disjoint from $Q_1 \cup Q_2$.
	\item $x_1, y_j \in V(Q_1)$; $x_2 \in V(Q_2)$; $x_i \in V(R_2)$, $i=2,\ldots,j$; $y_i \in V(R_2)$, $i=1,\ldots,j-1$.
	\item $y_0<y_1<x_3<y_2<x_4<\ldots <x_j<y_{j-1}<x_0$ in $R_2(y_0,x_0)$; $y_0<x_2<y_1$ in $R_2(y_0,x_0)$ or $x_2<y_0<y_1$ in $Q_2(v,u)$.
\end{enumerate}

Below, we assume that $y_0<x_2<y_1$ in $R_2(y_0,x_0)$ (the case when $x_2<y_0<y_1$ in $Q_2(v,u)$ is almost exactly the same; equation (\ref{eq(12)}) must be adjusted in the case $i=2$).

Set
\begin{equation}\label{eq(7)}
	C_i:=R_2(y_i,x_i)P_i(x_i,y_i), \quad (i=2,\ldots,j-1)
\end{equation}
Then $C_i$ is an $F$--alternating cycle.

Let $\overrightarrow{G^*}$ be a fixed even $F$--orientation of $G^*$ with associated even function $\omega$. Since $\omega$ is even, $\omega(C_i)=0$. Hence, from (\ref{eq(7)}),

\begin{equation}\label{eq(8)}
	\sum_{i=2}^{j-1}\omega(R_2(y_i,x_i))+\sum_{i=2}^{j-1} \omega(P_i(x_i,y_i)) \equiv 0 \, (\mbox{mod} \, 2) \,.
\end{equation}

\noindent {\sc Case $(a)$:} $x_1,y_j \in V(R_1)$.

Set
\begin{equation*}
	\begin{split}
		C_0:=Q_1(y_j,x_1)P_1(x_1,y_1)R_2(y_1,x_2)P_2(x_2,y_2)R_2(y_2,x_3) \ldots R_2(y_{j-2},x_{j-1}) & \\
		P_{j-1}(x_{j-1},y_{j-1})R_2(y_{j-1}x_j)P_j(x_j,x_j) \,.
	\end{split}
\end{equation*}
Then $C_0$ is an $F$--alternating cycle and $\omega(C_0)=0$. Hence,
\begin{equation}\label{eq(9)}
	\sum_{i=1}^{j}\omega(P_i(x_i,y_i))+\sum_{i=1}^{j-1} \omega(R_2(y_i,x_{i+1}))+ \omega(Q_1(y_j,x_1)) \equiv 0 \, (\mbox{mod} \, 2) \,.
\end{equation}
Also (see Lemma \ref{claim6.1} and its proof) because of the choice of $Q_1$, $Q_2$, $R_1$, $R_2$, $x_0$, $y_0$ and the maximality condition of Lemma \ref{claim6.1}(iv) (see Remark \ref{comment6.4} below), $\omega(C_i^*)=0$, $i=1,2$ where
\begin{equation*}
	\begin{split}
		C_1^*:=Q_1(x_1,y_0)R_2(y_0,y_1)P_1(y_1,x_1) & \\
		C_2^*:=Q_1(x_0,y_j)P_j(y_j,x_j)R_2(x_j,x_0) \,.
	\end{split}
\end{equation*}
Hence
\begin{equation}\label{eq(10)}
	\omega(Q_1(x_1,y_0))+ \omega(R_2(y_0,y_1))+\omega(P_1(y_1,x_1)) \equiv 0 \, (\mbox{mod} \, 2) \,,
\end{equation}
and
\begin{equation}\label{eq(11)}
	\omega(Q_1(x_0,y_j))+ \omega(P_j(y_j,x_j))+\omega(R_2(x_j,x_0)) \equiv 0 \, (\mbox{mod} \, 2) \,,
\end{equation}
Adding (\ref{eq(8)}), (\ref{eq(9)}), (\ref{eq(10)}) and (\ref{eq(11)}), we obtain:
\begin{equation}\label{eq(12)}
	\begin{split}
		(\sum_{i=1}^{j-1}\omega(R_2(y_i,x_{i+1}))+\sum_{i=2}^{j-1} \omega(R_2(y_i,x_i)))+ \omega(Q_1(y_j,x_1))+ & \\
		+\omega(Q_1(x_1,y_0))+ \omega(Q_1(x_0,y_j))+\omega(R_2(y_0,y_1))+\omega(R_2(x_j,x_0)) \equiv 0 \, (\mbox{mod} \, 2) \,.
	\end{split}
\end{equation}
Since $R_2(y_i,x_i)=R_2(y_i,x_{i+1})R_2(x_{i+1},x_i), $(i=2,\ldots, j-1), from (\ref{eq(12)}):
\begin{equation*}
	\begin{split}
		\omega(R_2(x_j,y_1))+ (\omega(Q_1(x_0,y_j)))+ \omega(Q_1(y_j,x_1))+\omega(Q_1(x_1,y_0))) & \\
		+(\omega(R_2(y_1,y_0))+1)+(\omega(R_2(x_0,x_j))+1) \equiv 0 \, (\mbox{mod} \, 2) \,.
	\end{split}
\end{equation*}
i.e.
\begin{equation}\label{eq(13)}
	\omega(R_1)+ \omega(R_2) \equiv 0 \, (\mbox{mod} \, 2) \,,
\end{equation}
which contradicts Lemma \ref{claim6.1}(iii).

\noindent{\sc Case $(b)$:} $x_1 \in V(R_1)$, $y_j \in V(Q_1(u,x_0))$.

The only difference from Case (a) is that now $C_2^*$ is an $F$--alternating cycle and hence $\omega(C_2^*)=0$, simply because $\omega$ is an even function.

\noindent{\sc Case $(c)$:} $x_1 \in V(Q_1(y_0,v))$, $y_j \in V(R_1)$.

This is the same as Case (b) up to a relabelling.

\noindent{\sc Case $(d)$:} $x_1 \in V(Q_1(y_0,v))$, $y_j \in V(Q_1(u,x_0))$.

This is the same as Case (a) except that now $\omega(C_i^*)=0$, $i=1,2$, simply since $\omega$ is an even function.\qed

\begin{remark}\label{comment6.4}
	Note that $\omega(C_i^*)=0$, $i=1,2$, by the maximality condition in Lemma \ref{claim6.1}(iv). For instance, consider the cycle $C_1^*$ and new paths $Q_1^*:=Q_1(u,v)$ and $Q_2^*:=Q_1(u,x_1)P_1(x_1,y_1)R_2(y_1,y_0)Q_1(y_0,v)$ with $R_1^*:=Q_1(x_1,y_0)$, $R_2^*:=P_1(x_1,y_1)R_2(y_1,y_0)$. By maximality $\omega(R_1^*)=\omega(R_2^*)$ i.e. $\omega(Q_1(x_1,y_0))=\omega(P_1(x_1,y_1))+\omega(R_2(y_1,y_0))$ (mod $2$). Since for odd length paths $P(u,v)$, $\omega(P(u,v))+\omega(P(u,v)) \equiv 1$ (mod $2$), we have $\omega(Q_1(x_1,y_0))+\omega(R_2(y_0,y_1))+\omega(P_1(y_1,x_1)) \equiv 0$ (mod $2$).
\end{remark}

\section{Proof of Theorem \ref{theorem3.5MAIN}}\label{section7}

%\Prf
Let $G_0$ be a minimal counterexample to Theorem \ref{theorem3.5MAIN}. From Theorem \ref{theorem5.10} and Lemma \ref{claim6.3}, $G^*=G_0-R$ where $R=\{e_1,e_2\}$ is a removable doubleton and $G^*$ is bipartite. Also $F$ is a fixed $1$--factor of $G_0$ such that $R\cap F=\emptyset$ and such that $G_0$ has no even $F$--orientations.

Let $\overrightarrow{G^*}$ be the canonical even $F$--orientation of $G^*$ with associated even function $\omega$ (cf. Definition \ref{defin2.9}). Assume that there do not exist cycles $C_1$ and $C_2$, relative to which $e_1$ and $e_2$ are skew. Let $e_1=(x_1,x_2)$ and $e_2=(y_1,y_2)$, $x_i \in X$ and $y_i \in Y$ $(i=1,2)$ and $(X,Y)$ be a bipartition of $G^*$. W.l.o.g., any cycle $C$ containing $e_1$ and $e_2$ is of the form
\setcounter{equation}{0}
\begin{equation}\label{eq(7.1)}
	C := (x_1,x_2,P_1(x_2,y_2), y_, y_1, P_2(y_1,x_1)).
\end{equation}
Since $\overrightarrow{G^*}$ is canonical, $\omega(P_1)=1$ and $\omega(P_2)=0$. Now define an $F$--alternating function $\omega_0$ on $G_0$ as follows:
\begin{enumerate}[(i)]
	\item if $(x,y) \in E(G_0^*)$, $\omega_0(x,y)=\omega(x,y)$;
	\item $\omega_0(x_1,x_2)=0$, $\omega_0(y_2,y_1)=1$.
\end{enumerate}
Then $\omega_0$ extends $\omega$ which itself is even. Hence, if $C$ is any cycle such that $R \cap E(C) = \emptyset$ then $\omega_0(C)=0$. If $R \cap E(C) \ne \emptyset$ then $R \subseteq E(C)$ and $C$ has the form of (\ref{eq(7.1)}). Then
\begin{equation*}
	\omega_0(C):=\omega_0(x_1,x_2)+\omega_0(P_1)+\omega_0(y_2,y_1)+\omega_0(P_2) \equiv 0 \,, (\mbox{mod} \, 2).
\end{equation*}
Hence, $\omega_0(C)=0$ for all $F$--alternating cycles $C$. Thus $G_0$ has an even $F$--orientation which is not true. Hence $G_0$ does have cycles $C_1$ and $C_2$ relative to which $e_1$ and $e_2$ are skew. Hence $G_0$ has a central subgraph $H$ $(H=G_0)$ such that $F$ is a $1$--factor of $H$ and $H$ is an even subdivision of a graph in $\cal{W}$. This contradicts the definition of $G_0$.\qed

\section{Proof of Theorem \ref{Thm2.4}}\label{section8}
%\section{Preliminaries to Theorem \ref{Thm2.4}}\label{section8}

%\textcolor{red}{Stare attenti a cambiare la conclusione delle dimostrazioni perché dicono ancora $F$-even}

In this section we introduce important tools, which will be useful in the proofs of Theorem \ref{Thm2.4}, which can be found at the end of the section.
%
%\begin{defin}\label{Weights} (Weights)
%
%\noindent Let $\overrightarrow{G}$ be an $F$-orientation of the graph $G$. Let $w$ be an additive $(0,1)$--function defined on the directed edges of $E(\overrightarrow{G})$ as follows:
%
%Let $\overrightarrow{P} \equiv (u_1,u_2, \ldots, u_n)$ denote an orientation of the $F$--alternating path $P(u_1,u_n)$. The ``opposite orientation'' of $\overrightarrow{P}$ is denoted by $\overleftarrow{P}$. Now define a function $w^*$ as follows. Set
%
%    \begin{equation*}
	%        w^*(u_i,u_{i+1}) =
	%        \left\{
	%        \begin{array}{ccc}
		%        1 & \mbox{if} & \overrightarrow{(u_i,u_{i+1})}   \\[16pt]
		%        0 & \mbox{if} & \overleftarrow{(u_i,u_{i+1})},
		%        \end{array}
	%        \right. \, , \, \, 1 \le i \le n
	%    \end{equation*}
%
%\noindent and $w^*(P) \equiv \sum_{i=1}^{n-1} w^*(u_i,u_{i+1}).$
%Similarly if $C$ is the $F$--alternating cycle
%$C:=(u_1,u_2, \ldots, u_n, u_1) = (P,u_1)$
%set $w^*(C) := w^*(P) + w^*(u_n,u_1)$. Finally set $w \equiv w^* (mod \, 2)$.
%
%We shall say that $w$ is the \emph{weight of the orientation} $\overrightarrow{G}$. \qed
%
%\end{defin}

\begin{lemma}\label{Lem3.2}
	Let $G \in {\cal W}$ with a ${\cal W}$-factor $F$. Let $\overrightarrow{G}$ be an $F$-orientation of $G$ with weight function $w$.
	An $F$-alternating cycle $C=(u_1,u_2, \ldots, u_n,u_1)$ is evenly oriented if and only if, for every $1 \le i \le n-1$
	\begin{equation*}
		w(u_1,u_2, \ldots, u_i) \equiv w(u_i,u_{i+1}, \ldots, u_n, u_1) \, ( \, mod \, 2).
	\end{equation*}
	
	%Let $w$ be the weight functions of $\overrightarrow{G}$, $G \in {\cal W}$, where $\overrightarrow{G}$ is an $F$--orientation of $G$. Let
	%$C=(u_1,u_2, \ldots, u_n,u_1)$ be an $F$--alternating cycle. Then, for
	%$1 \le i \le n-1 \, \, (i \mbox{ modulo } \, n)$
	%
	%\begin{equation*}
	%    w(u_1,u_2, \ldots, u_i) \equiv w(u_i,u_{i+1}, \ldots, u_n, u_1) \, ( \, mod \, 2)
	%\end{equation*}
	%
	%if and only if $C$ is evenly oriented.
\end{lemma}

\Prf
Follows immediately from the definitions.
\qed

\begin{lemma}\label{Lem3.3}
	Let $G \in {\cal W}$ with a ${\cal W}$-factor $F$. Let $\overrightarrow{G}$ be an $F$-orientation of $G$ with weight function $w$.
	If $P$ and $Q$ are $F$--alternating paths of odd length:
	%Let $\overrightarrow{G}$ be an even $F$--orientation of $G$, $G \in {\cal W}$, with weight function $w$.

	\begin{equation*}
		\begin{split}
			P := & (u_1,u_2, \ldots, u_k)  \\
			Q := & (v_1,v_2, \ldots, v_l)
		\end{split}
	\end{equation*}
	
	\noindent where $u_1=v_1,$ $u_k=v_l$ and $(u_1,u_2) \in F$, then $w(P) \equiv w(Q)$.
\end{lemma}

\Prf
Notice that since $(u_1,u_2) \in F$ and $k$ and $l$ are both even, $E(P \cup Q) \subset E(G) - R$ ($R$ as in Definition \ref{genwagner}). From Lemma \ref{Lem3.2}, if $E(P) \cap E(Q) = \emptyset$ the result immediately follows. So now assume that this is not the case.

Assume that the result is false.
Choose $P$ and $Q$ such that $w(P) \not\equiv w(Q)$ %$w(\overrightarrow{P}) \equiv w(\overrightarrow{Q})$
and such that $|V(P) \cap V(Q)|$ is minimal. Choose $j$ as small as possible such that $u_j \in V(P) \cap V(Q)$ ($j \ge 1$). Then, form above, $j < k$. And choose $i$ so that $i \le j-1$, $P(u_i,u_j) \subset Q$ and $i$ is as small as possible. From Lemma \ref{Lem3.2}
\begin{equation}\label{Eq1}
	w(P(u_1,u_i)) = w(\overleftarrow{Q}(v_n,v_1))
\end{equation}
\noindent where $v_n=u_i$.

Now replace $u_1$ by $u_i (=v_n)$ in the above argument and replace the paths $P$ and $Q$ by $P(u_i,u_k)$ and $Q(v_m,v_k)$ respectively, using (\ref{Eq1}) and in minimality we obtain a contradiction.
\qed

\begin{lemma}\label{Lem3.4}
	Let $G \in {\cal W}$ with a ${\cal W}$-factor $F$, and $u,v \in V(G)$. Let $\overrightarrow{G}$ be an $F$--orientation of $G$.
	Then $G$ contains an $F$--alternating path $P^*(u,v)$ of odd length with first edge not in $F$.
	\noindent
\end{lemma}

%\begin{equation*}
%\begin{split}
%    P(u,v) := & (u_1,u_2, \ldots, u_k) \, ; \, (u_1,u_2) \in E(G) \setminus F \\
%              & u_1 = u, u_k=v~and~k \ge 2,~k~even.
%\end{split}
%\end{equation*}

\Prf
Choose $P:=P(u,v):=(u_1,u_2, \ldots, u_k)$, $u_1 = u, u_k=v~and~k \ge 2,~k~even$ with $(u_1,u_2) \in E(G) \setminus F$ and
such that $|E(P) \cap F|$ is as large as possible. Now choose $i$ as large as possible such that $P(u_1,u_i)$ is $F$--alternating.
If $i=k$, we are done. Otherwise, set $e_0=(u_i,u_{i+1})$. Since $G$ is $1$-extendible, it is $2$--connected, so there exists at least
one cycle in $G$ containing $e_0$. If $e_0 \in R$ ($R$ as in Definition \ref{genwagner}), choose $C \in \{C_1,C_2\}$. Otherwise, let $C \in G-R$ be a cycle containing $e_0$, which is $F$-alternating for the bipartiteness of $G-R$ .
In both cases, $C$ is an $F$--alternating cycle containing $e_0$.
Then $P \cup C$ contains a path $P':=P'(u,v)$ such that
$|E(P') \cap F| > |E(P) \cap F|$ which is a contradiction to the choice of $P$. This implies that $i$ must be equal to $k$
and we can set $P^*(u,v)=P(u,v)$.
\qed

\begin{lemma}\label{Lem3.5}
	Let $G \in {\cal W}$ with a ${\cal W}$-factor $F$, and $u,v \in V(G)$. Let $\overrightarrow{G}$ be an $F$--orientation of $G$.
	Let $P(u,v)$ be an $F$--alternating path in $\overrightarrow{G}$ with first and last edges in $F$. Then if $Q(u,v)$ is any $F$--alternating path
	
	\begin{equation*}
		w(P) = w(\overleftarrow{Q}).
	\end{equation*}
	
\end{lemma}

\Prf
This follows from Lemmas \ref{Lem3.3} and \ref{Lem3.4}.
\qed

\begin{figure}[h]
	\begin{center}
		\includegraphics[width=9cm]{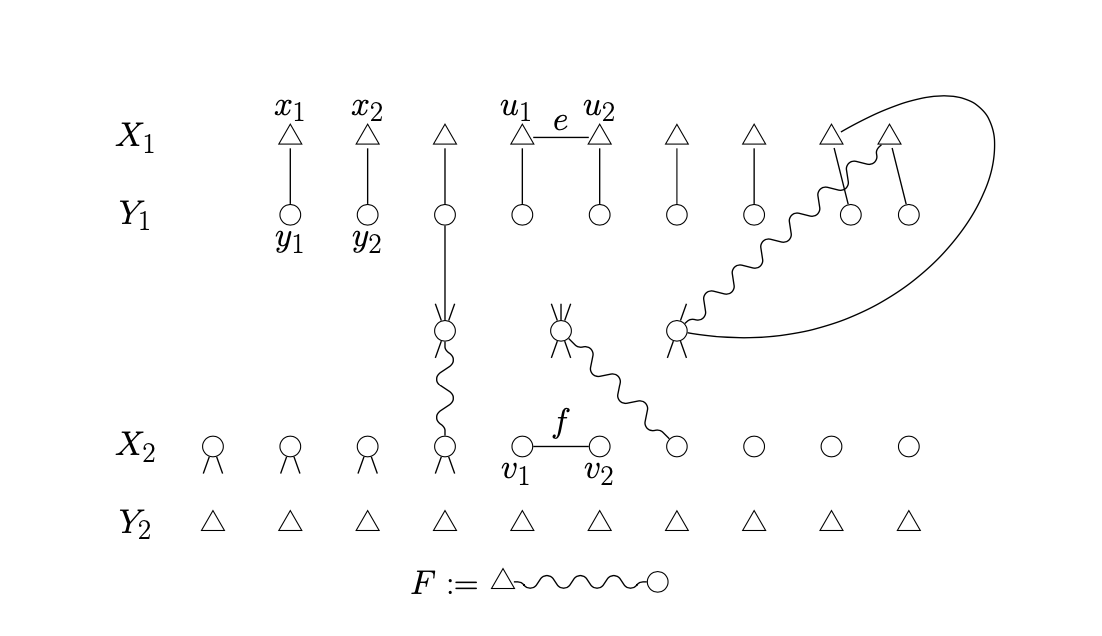}\\
	\end{center}
	\caption{Illustration of Notation \ref{Not5.1}}\label{Fig5.2}
\end{figure}

\begin{notation}\label{Not5.1}
	Suppose that $G \in {\cal {W}}$ (see Figure \ref{Fig5.2}, and set $\ell := \kappa (G)$.
	Let $S = \{s_1,s_2, \ldots, s_{\ell}\}$ be a separating set. Let $G \setminus S:= G_1 \cup G_2$ where $e \in E(G_1)$ and $f \in E(G_2)$. Suppose that
	$V(G_1) \setminus \{e\}:=X_1 \cup Y_1$ and
	$V(G_2) \setminus \{f\}:=X_2 \cup Y_2$ where
	$X=X_1 \cup X_2$ and $Y=Y_1 \cup Y_2$.
	
	Set $F_i:= F \cap E(S, X_i \cup Y_i)$ $(i=1,2)$ and finally set
	$F_{i_1}:= F_i \cap E(S, X_i)$ and $F_{i_2}:= F_i \cap E(S, Y_i)$ $(i=1,2)$.
	
	Set $|X_i|=\xi_i$, $|Y_i|=\tau_i$; $k_{ij} := |E(s_i,X_j \cup Y_j)|$
	$(i=1,2,\ldots,\ell; j=1,2)$.
\end{notation}

\Prf \textbf{of Theorem \ref{Thm2.4}(i)}

\begin{figure}[h]
	\begin{center}
		\includegraphics[width=6cm]{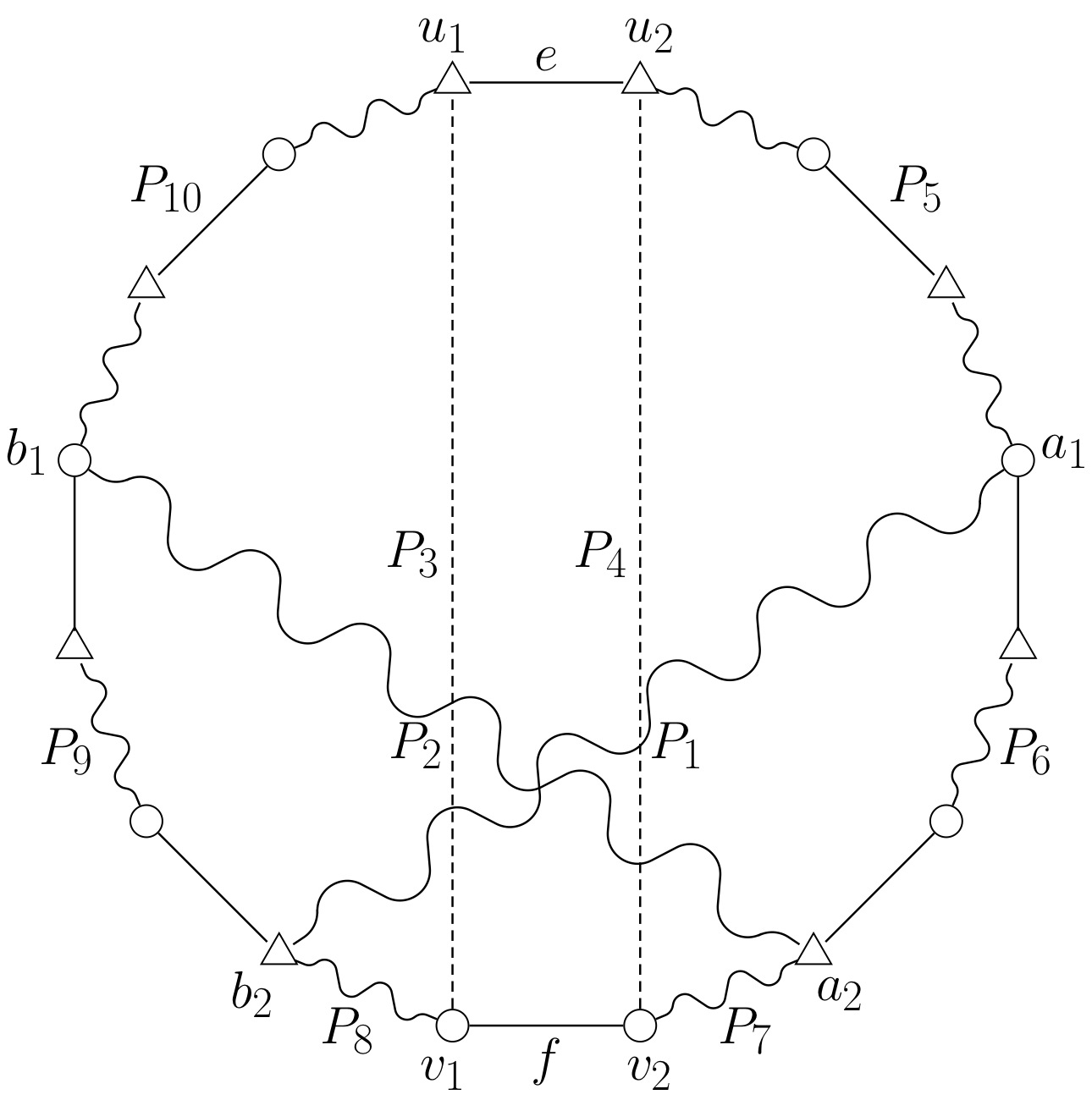}\\
	\end{center}
	\caption{Even subdivision of $K_4$}\label{Fig5.3}
\end{figure}

Suppose that $G \in {\cal{W}}$ and $\kappa (G) = 4$. By Lemma \ref{lemma5.2}, $G$ contains an $F$--central subgraph $H \in {\cal W}(\le 3)$ which is isomorphic to an even subdivision of $K_4$ (cfr. Figure \ref{Fig5.3}).

In Figure \ref{Fig5.3} $P_1$ and $P_2$ ($P_1 \cap P_2 = \emptyset$) denote paths $P_1 := P_1(a_1,b_2)$, $P_2 := P_2(b_1,a_2)$, $a_1,b_1 \in X$ and
$a_2,b_2 \in Y$.
We set $e:=(u_1,u_2)$ and $f:=(v_1,v_2)$. The skew cycles $C_1$ and $C_2$ are

$$C_1:=(u_1,u_2, \ldots,a_1, \ldots, a_2, \ldots, v_2,v_1, \ldots, b_2,b_1, \ldots, u_1)$$
$$C_2:=(u_1,u_2, \ldots,a_1,P_1,b_2, \ldots, v_1,v_2, \ldots, a_2,P_2,b_1, \ldots, u_1).$$

Since $\kappa (G)=4$ there exist disjoint paths $P_3^*:=P_3^*(u_1,v_1)$ and $P_4^*:=P_4^*(u_2,v_2)$ (or we can relabel say $u_1$ and $v_1$). Now from Lemma \ref{Lem3.4} there exist $F$--alternating paths $P_3:=P_3(u_1,v_1)$ and $P_4:=P_4(u_2,v_2)$.
%In Figure \ref{Fig5.3} we draw these paths (although in reality they may well not be disjoint).
Now let us suppose that $G$ has an $F$--even orientation $\overrightarrow{G}$ with weight function $w$. Set (again see Figure \ref{Fig5.3})

\begin{equation*}
	\begin{array}{ccc}
		P_5:=P_5(u_2,\ldots,a_1) & P_6:=P_6(a_1,\ldots,a_2) & P_7:=P_7(a_2,v_2) \\
		P_8:=P_8(v_1,b_2) & P_9:=P_9(b_2,\ldots,b_1) & P_{10}:=P_{10}(b_1,\ldots,u_1)
	\end{array}
\end{equation*}

Now using Lemmas \ref{Lem3.3} and \ref{Lem3.5} and setting $w_i:=w(P_i)$ ($i=1,2, \ldots, 10$)

\begin{equation*}
	\begin{split}
		w(u_1,u_2) + w_5 + w_6 + w_7 + w(v_2,v_1) + w_8 + w_9 + w_{10} \equiv 0 \quad (mod \quad 2)\\
		w(u_1,u_2) + w_5 + w_1 + \overline{w_8} + w(v_1,v_2) + \overline{w_7} + w_2 + w_{10} \equiv 0 \quad (mod \quad 2)\\
		w_3 + w_8 + w_9 + w_{10} \equiv 0 \quad (mod \quad 2)\\
		w_5 + w_6 + w_7 + \overline{w_4} \equiv 0 \quad (mod \quad 2)\\
		w_3 + w_8 + \overline{w_1} + \overline{w_5} + w_4 + \overline{w_7} + w_2 + w_{10} \equiv 0 \quad (mod \quad 2)
	\end{split}
\end{equation*}

where $\overline{w_i}:= 1 + w_i$ ($i=1,2, \ldots, 10$)

Adding this equation:

\begin{equation*}
	w_5 + w(v_1,v_2) + w_8 + w_1 + w_4 \equiv 0 \quad (mod \quad 2)
\end{equation*}

which is a contradiction. \qed

\begin{lemma}\label{Lem5.4}
	Let $G \in {\cal W}$ and $G$ be regular of degree $k$ ($k \ge 3$). Suppose that $\kappa (G) = 2$. Then $G$ is not $\cal{F}$--even.
\end{lemma}

\Prf
Suppose that $G \in {\cal W}$, $\kappa (G) = 2$ and $G$ $k$--regular ($k \ge 3$). Suppose that $G$ is $\cal{F}$--even. Set $S:=\{s_1,s_2\}$ where $S$ is a separating set. There are (several cases) to consider.

{\sc Case $1$} ($|F_{11}| + |F_{12}| \equiv 0 \, (mod \, 2)$)
In these case $G \notin {\cal W}$ since the skewness condition of ${\cal W}$ is contradicted.

{\sc Case $2$} ($|F_{11}| = 1,  |F_{12}| = 0$)
In this case $\xi_1 = \tau_1 + 1$ and $k\xi_1 - k_{11} - 2 = k\tau_1 - k_{12}$. Hence $k = 2 + k_{11} + k_{12}$ and $k_{11} = k - 1$, $k_{12} = 1$. Set $E_0 = \{(s_1,y_1),(s_2,x_2)\}$ where $y_1 \in Y_1$, $x_2 \in X_2$, $s_1 \in X$, $s_2 \in Y$. Then $E_0$ is an edge--cut. Let $G_1 \setminus E_0 = H_1 \dot\cup H_2$ and set $H_1^*:=H_1 + (s_1,x_2)$ and $H_2^*:=H_2 + (y_1,s_2)$. Clearly $H_i^* \in {\cal W}$ for some $i \in \{1,2\}$. Furthermore we can choose $H_i^*$ to be as small as possible. In particular we can choose $H_i^*$ so that $\kappa (H_i^*) \ge 3$ and $H_i^*$ is $k$--regular ($k \ge 3$). So by {\sc Case $1$}, $H_i^*$ is not $\cal{F}$--even which implies $G$ is not $\cal{F}$--even which is not true.

{\sc Case $3$} ($|F_{11}| = 2$)
In this case $\xi_1 = \tau_1 + 2$ and $\xi_1k - k_{11} - k_{12} - 2 = \tau_1k.$ Hence $2k = k_{11} + k_{12}+2$. Hence $k_{11}=k_{12}=k-1$. Set $E_0 = \{(s_1,y_1),(s_2,y_2)\}$ where $y_i \in Y_2$, $s_i \in X$ ($i=1,2$). This $E_0$ is an edge--cut and we repeat the argument of {\sc Case $2$}. \qed

%The reader might find the following proof of Theorem \ref{theorem3.5MAIN}(iii) easier to digest after reading Example \ref{Ex5.8} at the end of this paper.

\Prf \textbf{of Theorem \ref{Thm2.4}(ii)}
Let $G \in {\cal W}$ be regular of degree $k$ ($\geq 3$). Assume that $G$ is $\cal{F}$--even. From Theorem \ref{Thm2.4}(i) and Lemma \ref{Lem5.4}, $\kappa (G) = 3.$

We use the terminology of the introduction and of Notation \ref{Not5.1}. Thus $S = \{s_1,s_2,s_3\}$ is a separating set. By symmetry we may assume that $|S \cap X| \ge 2$. We now prove, with this assumption, that

\setcounter{equation}{0}
\begin{equation}\label{Eq1Prf2.6}
	|S \cap X| \equiv |F_{11}| + |F_{12}| \quad (mod  \quad 2)
\end{equation}

Set $k_{ij}:=|E(s_i,X_, \cup Y_j|$ ($i=1,2,3; \, j=1,2$). There are two cases

{\sc Case $1$} ($|S \cap X| = 3$)
Since $|S \cap (V(C_1) \cup V(C_2))| \ge 2$, $|F_{12}| \le 1$. Suppose that $|F_{12}|=0.$ Then $|F_{22}| = 3$. Hence $\tau_2 = \xi_2 + 3$ and it follows that $\tau_2k - 2 - (k_{21} + k_{22} + k_{23}) = \xi_2k$, $3k = 2 + (k_{21} + k_{22} + k_{23}) \le 3k - 1$, $|F_{12}|=1$ and, by definition, $|F_{11}|=0$. Therefore (\ref{Eq1Prf2.6}) is satisfied.

{\sc Case $2$} ($|S \cap X| = 2$)
Suppose that $s_1,s_2 \in X$ and $s_3 \in Y$. As in the previous case $|F_{12}| \le 1$. By definition $|F_{11}| \le 1$.

Suppose that $|F_{12}| = 0$. Assume that $|F_{11}| = 1$. Then $\tau_2 = \xi_2 + 2$. Hence $\xi_2k - k_{23} = \tau_2k - 2 - k_{21} - k_{22}$. Hence $2k = 2 + k_{21} + k_{22} - k_{23} \le 2k -1$. Hence if $|F_{12}| = 0$, $|F_{11}|=0$ and (\ref{Eq1Prf2.6}) is satisfied.

Finally suppose that $|F_{12}| = 1$. Assume that $|F_{11}|=0$. Then this is impossible since $|S \cap (V(C_1) \cup V(C_2))| \ge 2$. Hence $|F_{11}|=1$.

Hence in all cases (\ref{Eq1Prf2.6}) is satisfied.

\

Now we assume that $|S \cap X| = 2$ and $|F_{11}|=|F_{12}|=1$ (see (\ref{Eq1Prf2.6}) above. We follow the method of proof of Theorem \ref{Thm2.4}(i). By Lemma \ref{lemma5.2}, $G$ contains an $F$--central subgraph $H \in {\cal W}( \le 3)$ which is isomorphic to an even subdivision of $K_4$ (See Figure \ref{Fig5.2} in the proof of Theorem \ref{Thm2.4}(i). %We follow the argument and notation of Theorem \ref{Thm2.5} up to the point when in this proof we must adjust the assumption of connectivity to our assumption of $\kappa (G) = 3.$

\begin{figure}[h]
	\begin{center}
		\includegraphics[width=12cm]{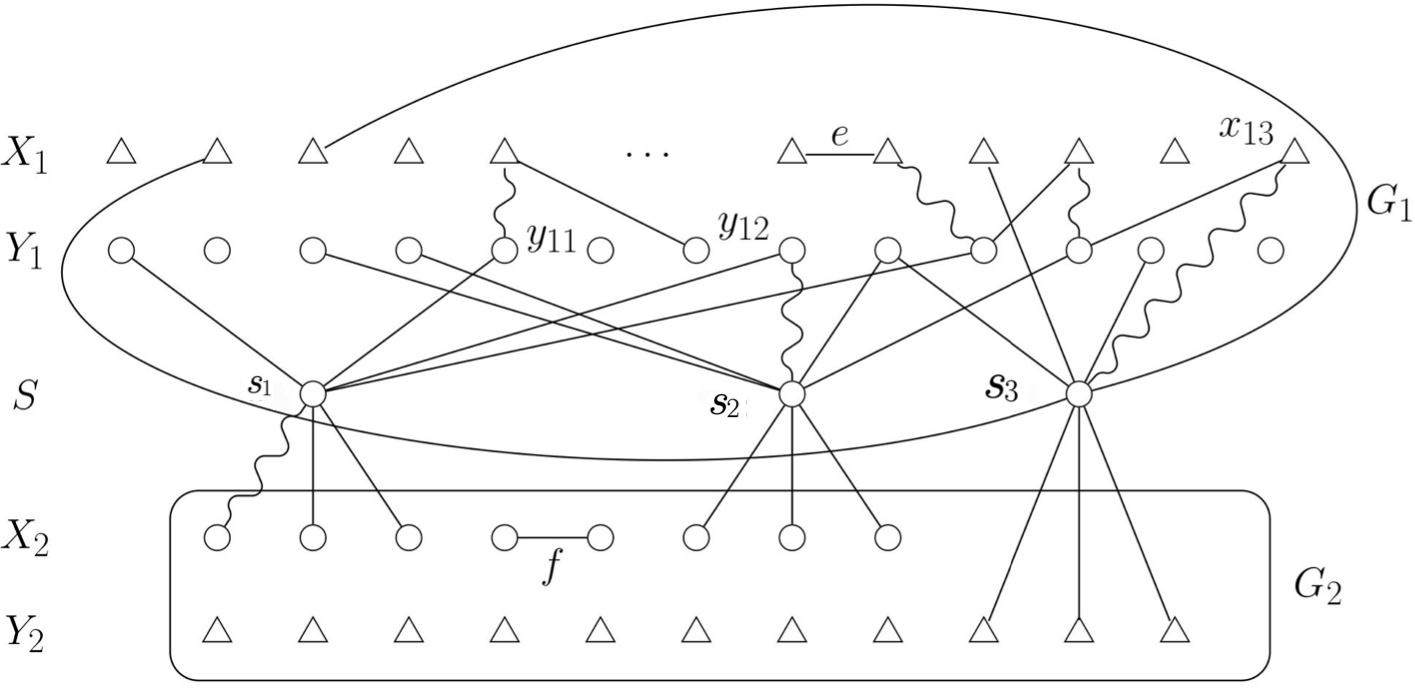}\\
	\end{center}
	\caption{Illustration for the proof of Theorem \ref{Thm2.4}(ii)}\label{Fig5.5}
\end{figure}

We may assume (see Figure \ref{Fig5.5} without loss of generality that there exist edges $s_1,y_{11} \in E(G) \setminus F$ and $(s_2,x_{13}) \in F_1$, $y_{11} \in Y_1$, $x_{13} \in X_1$, belonging to $E(C_1) \cup E(C_2).$ Also there exists an edge $(s_2,y_{12}) \in F$, $y_{12} \in Y_1.$

Suppose that $e$ and $y_{12}$ are in different components of $G_1$. If $y_{12}$ and $x_{13}$ are in the same component then since $x_{13} \in V(C_1) \cup V(C_2)$, $y_{12}$ and $e$ are in the same component which is a contradiction.

So now assume that $x_{13}$ and $y_{12}$ belong to different components of $G_1$. Consider the component of $G_1$ containing $x_{13}$ (which must clearly contain $e$) having sides $X_{12} \subseteq X_1$, $Y_{12} \subseteq Y_1$,
$\xi_{12} :=|X_{12}|$ and $\tau_{12} := |Y_{12}|$. Then $\xi_{12}=\tau_{12} + 1$ and $\xi_{12}k - 2 - k^*_{23} = \tau_{12k} - k^*_{21} - k^*_{22}$ where $k^*_{ji} := |E(s_i, X_{12} \cup Y_{12}|.$ Hence

$$k^*_{21} + k^*_{22} = 2 + k^*_{23} - k \le 1.$$

Hence there exists a proper subset of $S$ which separates $K$ form $G \setminus K$ which is a contradiction. Hence $x_{13}$ and $y_{12}$ belong to the same component of $G_1$ which also includes $e$.

%Now again returning to the method of proof of Theorem \ref{theorem3.5MAIN}(iii) (\textbf{quale sarebbe?})
Without loss of generality, there exist paths $P^*_{13}:=P^*_{13}(u_1,y_{12})$ and $P^*_{14}:=P^*_{14}(u_2,y_{12})$ in $G_1 \setminus \{e\}$ with both paths having their final edges in $E(G) \setminus F$ (note that the notation may be chosen so that $u_1$ and $u_2$ or $v_1$ and $v_2$ may be interchanged). Similarly in $G \setminus G_1$ there exist paths
$P^*_{23}:=P^*_{23}(s_2,v_1)$ and $P^*_{24}:=P^*_{24}(s_2,v_2).$ Finally set $P^*_3:=P^*_{13}P^*_{23}$ and $P^*_4:=P^*_{14}P^*_{24}$ and continue exactly as in Theorem \ref{Thm2.4}(i).

\begin{figure}[h]
	\begin{center}
		\includegraphics[width=12cm]{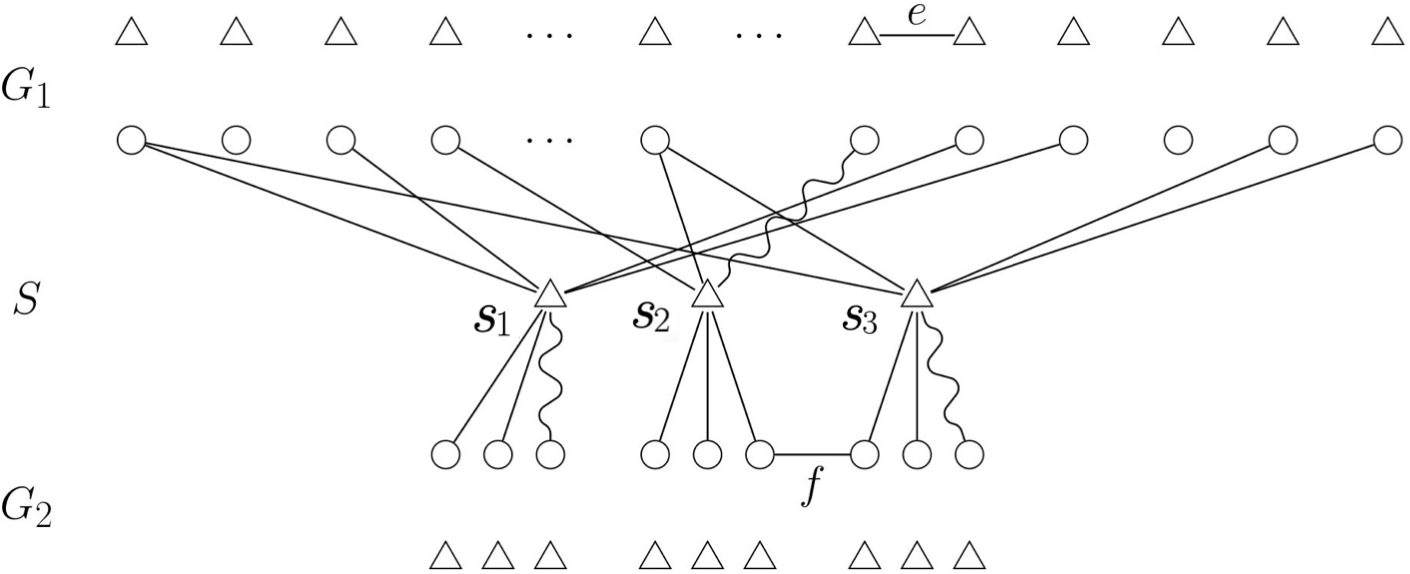}\\
	\end{center}
	\caption{Illustration for the proof of {\sc Case A} in Theorem \ref{Thm2.4}(ii)}\label{Fig5.6}
\end{figure}

There are now two other cases to consider. In fact these cases basically duplicate the first case:

{\sc Case A} ($|S \cap X|=3$, $|F_1|=1$, $|F_{11}|=0$)
Suppose that $(s_2,y_{12}) \in F_1$, $y_{12} \in Y_1$ (see Figure \ref{Fig5.6})

Suppose that $y_{12}$ and $e$ are in different components of $G_1$. Then $\tau_2=s_2$ and $\xi_2k - 2 = \tau_2k- k^*_{21}-k^*_{22}-k^*_{23}.$ Hence $k^*_{21}+k^*_{22}+k^*_{23}=2$. Hence the component $K$, say, containing $e$ in $G_1$ is separated from $G \setminus K$ by a proper subset of $S$ which is a contradiction.

{\sc Case B} ($|S \cap X|=2$, $|F_{11}|=|F_{12}|=0$)

\begin{figure}[h]
	\begin{center}
		\includegraphics[width=12cm]{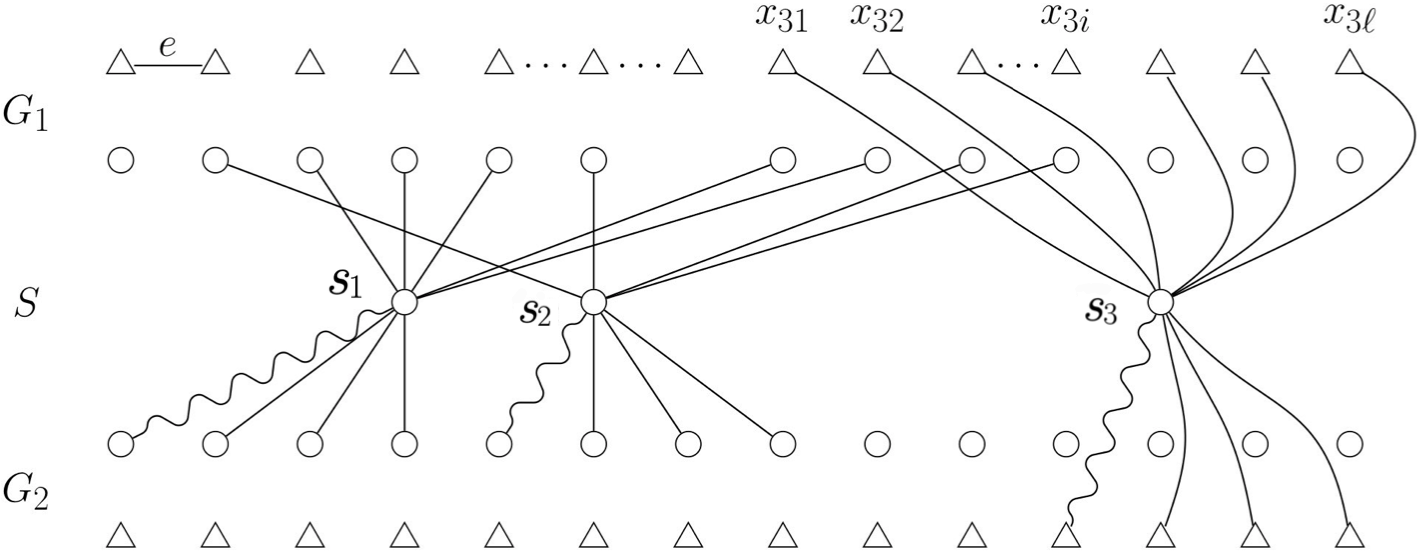}\\
	\end{center}
	\caption{Illustration for the proof of {\sc Case B} in Theorem \ref{Thm2.4}(ii)}\label{Fig5.7}
\end{figure}

Suppose that $(s_3,x_{3i}) \in E(G) \setminus F$, $x_{3i} \in X_1$, $i=1,2, \ldots, \ell$ (see Figure \ref{Fig5.7})

Now suppose that there is no component in $G_1$ containing both $e$ and $x_{3i}$ for any $i \in \{1,2 , \ldots, \ell\}$. Again this would imply that the component $K$ containing $e$ in $G_1$ is separated from $G \setminus K$ by a proper subset of $S$ which is the final contradiction.
\qed

\bigskip

\begin{example}\label{Ex5.8}
	We give a concrete example illustrating Theorem \ref{Thm2.4}(ii). The graph $G$ in Figure \ref{Fig5.9} is such that $G \in {\cal W}$, $\kappa(G) =3$ and $G$ is $4$--regular $G$ has a separating set $S=\{s_1,s_2,s_3\}$ where $s_1 \in X$, $s_2,s_3 \in Y.$

	\begin{figure}[h]
		\begin{center}
			\includegraphics[width=6cm]{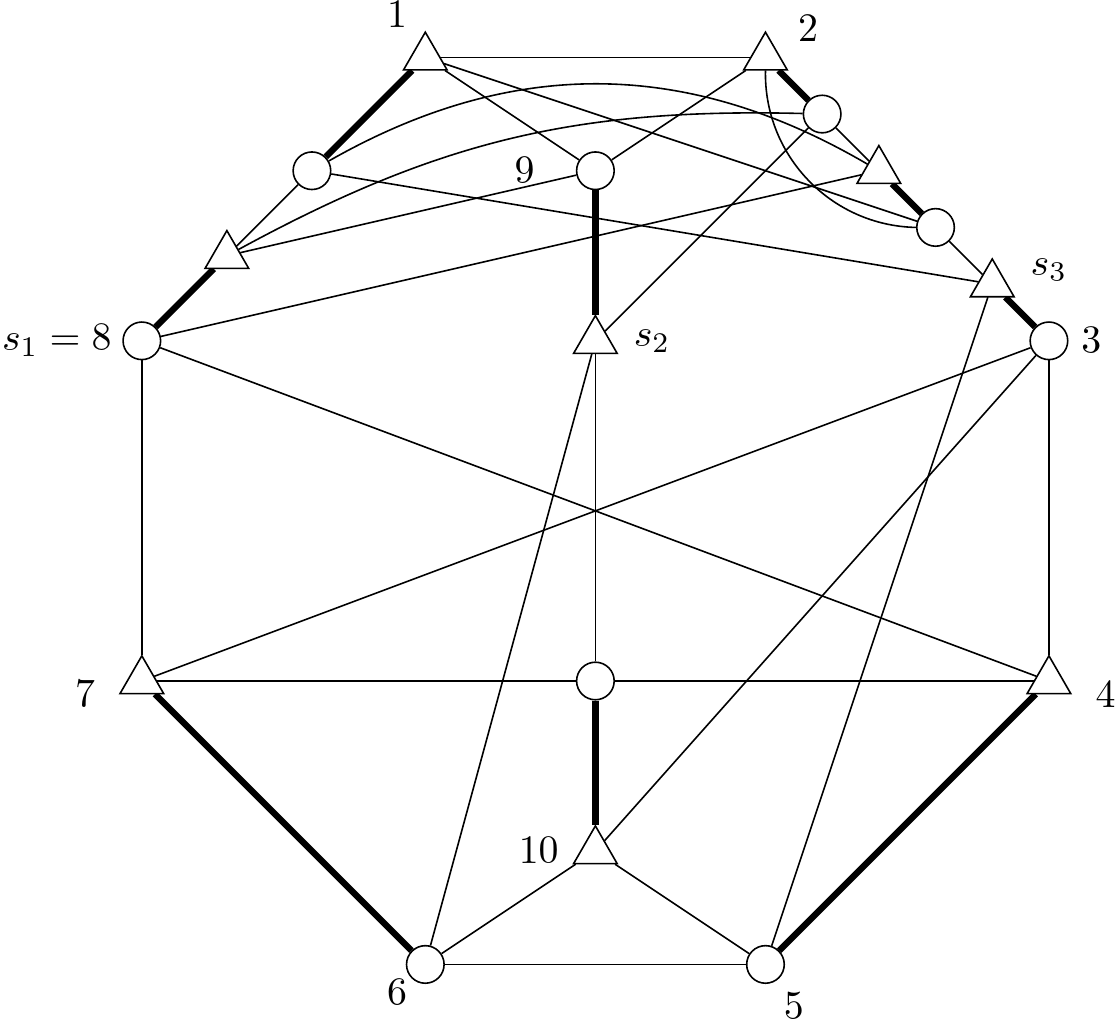}\\
		\end{center}
		\caption{$G$}\label{Fig5.9}
	\end{figure}
	
	In Figure \ref{Fig510} the graph $G_0$ is an $F$--central subgraph of $G$ and $\overrightarrow{G_0}$ is an $F$--orientation. We use the labelling of Figure \ref{Fig5.9} except now $s_1$ is relabelled $8$.
	
	\begin{figure}[h]
		\begin{center}
			\includegraphics[width=4cm]{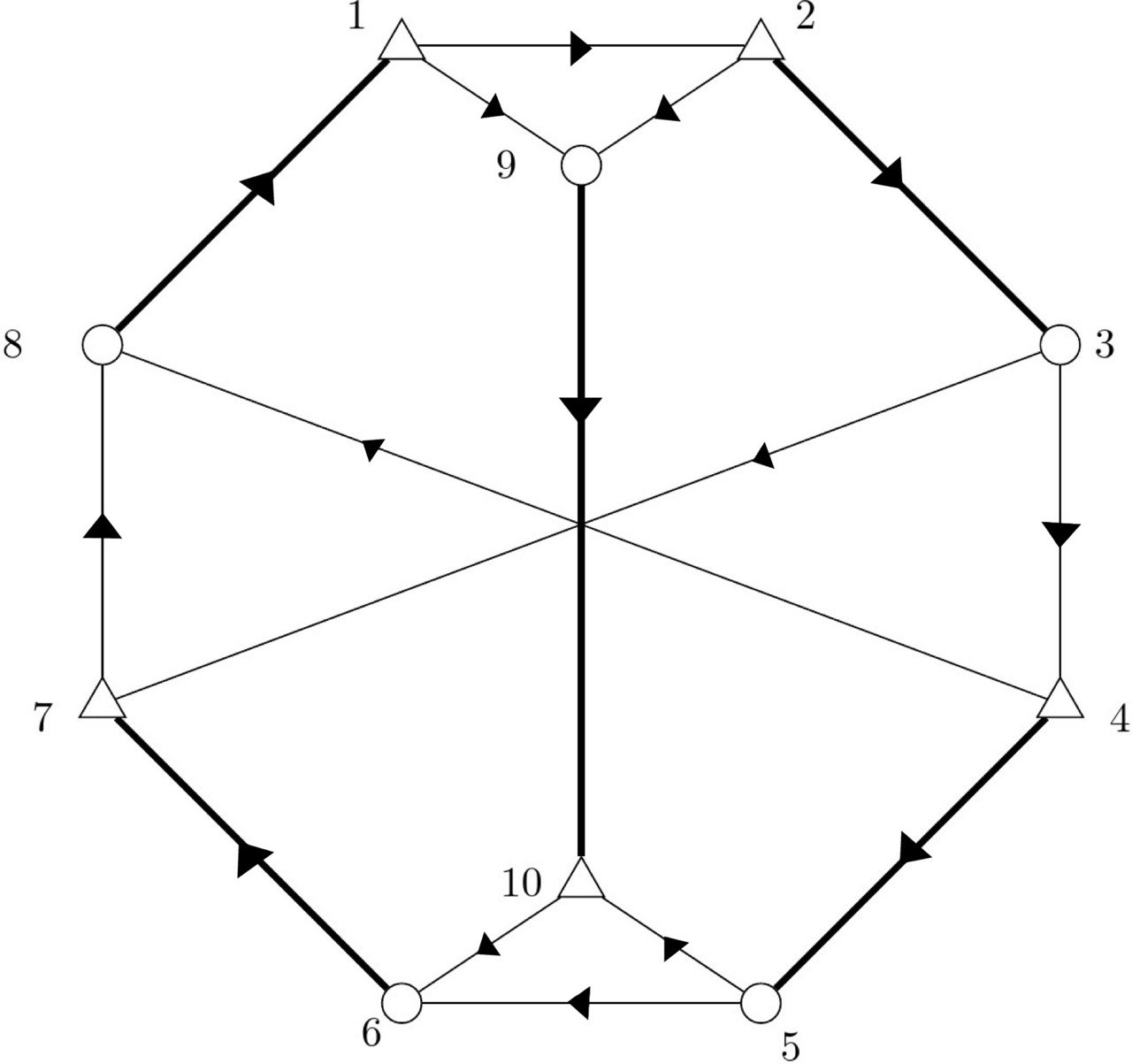}\\
		\end{center}
		\caption{$G_0$ = an $F$--central subgraph of $G$}\label{Fig510}
	\end{figure}

	$G_0$ has $F$--alternating cycles:
	
	\begin{equation*}
		\begin{array}{lll}
			C_1:=(1,2,3,4,5,6,7,8,1), & \quad & C_2:=(1,9,10,6,7,8,1), \\
			C_3:=(2,9,10,5,4,3,2), &  & C_4:=(1,9,10,5,4,8,1), \\
			C_5:=(2,9,10,6,7,3,2), & & C_6:=(1,2,3,7,6,5,4,8,1).
		\end{array}
	\end{equation*}

	Then $\{C_1,C_2, \ldots, C_6\}$ is a zero--sum set and in this set $C_6$ is the only evenly $F$--oriented cycle. This proves that $G_0$ is not $F$--even and hence $G$ is not $\cal{F}$--even.
	\qed
\end{example}

\section*{Acknowledgements}
Thanks to the group GNSAGA of INdAM for their funding and support.

\bibliography{bibtex4even}{}
\bibliographystyle{plain}

\end{document}